\documentclass[16pt]{article}
\usepackage{amsmath,amsfonts,amssymb}
\usepackage{verbatim}
\usepackage{latexsym}
\newtheorem{theorem}{Theorem}[section]
\newtheorem{lemma}[theorem]{Lemma}
\newtheorem{proposition}[theorem]{Proposition}
\newtheorem{definition}[theorem]{Definition}

\newcommand\RR{{\Bbb R}}
\newcommand\CC{{\Bbb C}}

\newcommand\NN{{\Bbb N}}
\newcommand\ZZ{{\Bbb Z}}

\newcommand\NT{{\Bbb N}}

\newcommand\supp{{\rm supp}}

\begin{document}
\title{Besov Spaces and Frames on Compact Manifolds }
\author{Daryl Geller\\
\footnotesize\texttt{{daryl@math.sunysb.edu}}\\
 Azita Mayeli \\
\footnotesize\texttt{{amayeli@math.sunysb.edu}}
}

\maketitle

\begin{abstract}
We show that one can characterize the Besov 
spaces on a smooth compact oriented Riemannian manifold, for the full range of indices, through a knowledge of the size of frame
coefficients, using the smooth, nearly tight frames we have constructed in \cite{gm2}. 

\end{abstract}

\section{Introduction}
In \cite{gm2}, we have constructed smooth, nearly tight frames on $({\bf M},g)$, a general smooth, compact oriented 
Riemannian manifold without boundary.
Our goal in this article is to show that one can characterize the (inhomogeneous) Besov 
spaces on ${\bf M}$, for the full range of indices, through a knowledge of the size of frame
coefficients, using the frames we have constructed.  (We hope to consider Triebel-Lizorkin spaces in a
forthcoming article.)  Our methods, in addition to using the results of \cite{gm2},
are largely adapted from those of Frazier and Jawerth \cite{FJ1},
who gave a similar characterization of Besov spaces on $\RR^n$.   However, as we shall explain below,
some new ideas are needed on manifolds.

Let us briefly review our construction of smooth, nearly tight frames on ${\bf M}$.  

Say $f^0 \in {\cal S}(\RR^+)$ (the space of restrictions to $\RR^+$ of functions in ${\mathcal S}(\RR)$).
Say $f^0 \not\equiv 0$, and let 
\[f(s) = sf^0(s).\]
One then has
the {\em Calder\'on formula}: if $c \in (0,\infty)$ is defined by
\[c = \int_0^{\infty} |f(t)|^2 \frac{dt}{t} = \int_0^{\infty} t|f^0(t)|^2 dt,\]
then for all $s > 0$,

\begin{equation}
\label{cald}
\int_0^{\infty} |f(ts)|^2 \frac{dt}{t} = c < \infty.
\end{equation}

Discretizing (\ref{cald}), if $a \neq 1$ is sufficiently close to $1$,
one obtains a special form of
{\em Daubechies' condition}:   
for all $s > 0$,
\begin{equation}
\label{daub}
0 < A_a \leq \sum_{j=-\infty}^{\infty} |f(a^{2j} s)|^2 \leq B_a < \infty,
\end{equation}
where  
\begin{equation}
\label{daubest}
A_a = \frac{c}{2|\log a|} \left(1 - O(|(a-1)^2 (\log|a-1|)|\right),\:\:\: B_a=
\frac{c}{2|\log a|}\left(1 + O(|(a-1)^2 (\log|a-1|)|)\right).
\end{equation}
((\ref{daubest}) was proved in \cite{gm1}, Lemma 
7.6)  In particular, $B_a/A_a$ converges nearly quadratically to $1$ as $a \rightarrow 1$.  For example,
Daubechies calculated that if $f(s) = se^{-s}$ and $a=2^{1/3}$, then $B_a/A_a = 1.0000$  to four 
significant digits.

Our general program is to construct (smooth, nearly tight) frames, and analogues of continuous wavelets, on 
much more general spaces, by replacing the positive number $s$ in (\ref{cald}) and
(\ref{daub}) by a positive self-adjoint operator $T$ on a Hilbert space ${\cal H}$.  
If $P$ is the projection onto the null space of $T$, by the spectral theorem 
we obtain the relations

\begin{equation}
\label{gelmay1}
\int_0^{\infty} |f|^2(tT) \frac{dt}{t} = c(I-P)
\end{equation}

and

\begin{equation}
\label{gelmay2}
A_a(I-P) \leq \sum_{j=-\infty}^{\infty} |f|^2(a^{2j} T) \leq
B_a(I-P). 
\end{equation}

(The integral in (\ref{gelmay1}) and the sum in (\ref{gelmay2}) converge strongly.
In (\ref{gelmay2}), $\sum_{j=-\infty}^{\infty} := 
\lim_{M, N \rightarrow \infty} \sum_{j=-M}^N$, taken in the strong operator 
topology.) 
(\ref{gelmay1}) and (\ref{gelmay2}) were justified in section 2 of our earlier article 
\cite{gmcw}. 

In \cite{gmcw} and \cite{gm2}, we looked at the
situation in which $T$ is the Laplace-Beltrami operator $\Delta$ on $L^2({\bf M})$, 
We constructed smooth, nearly tight frames
in this context.  Here $P$ is the projection onto the 
one-dimensional space of constant functions.  
We constructed continuous wavelets on ${\bf M}$ in \cite{gmcw}.

To see how frames can be obtained from (\ref{gelmay2}), suppose that, for any $t > 0$,
$K_t$ is the Schwartz kernel of $f(t^2T)$.  Thus, if $F \in L^2({\bf M})$,  

\begin{equation}
\label{schkerf}
[f(t^2T)F](x) = \int_{\bf M} F(y) K_t(x,y) d\mu(y), 
\end{equation}
here $\mu$ is the measure on ${\bf M}$ arising from integration with respect to the
volume form on ${\bf M}$.  
Say now that $\int_{\bf M} F = 0$, so that $F = (I-P)F$.
By (\ref{gelmay2}),
\begin{equation}
\label{aasumba0}
A_a \langle F,F \rangle  \leq  \langle \sum_j |f|^2(a^{2j} T)F, F \rangle  \leq B_a \langle F,F \rangle .
\end{equation}
Thus
\begin{equation}
\label{aasumba}
A_a \langle F,F \rangle  \leq \sum_j  \langle f(a^{2j} T)F, f(a^{2j} T)F \rangle  \leq B_a \langle F,F \rangle ,
\end{equation}
so that
\begin{equation}
\label{kerfaj}
A_a \langle F,F \rangle  \leq \sum_j \int|\int K_{a^j}(x,y) F(y)d\mu(y)|^2 d\mu(x) \leq B_a \langle F,F \rangle .
\end{equation}
Now, pick $b > 0$, and for each $j$, write ${\bf M}$
as a disjoint union of measurable sets $E_{j,k}$ with diameter at most $ba^j$.
Take $x_{j,k} \in E_{j,k}$.
It is then reasonable to expect that, for any $\epsilon > 0$, if 
$b$ is sufficiently small, and if $x_{j,k} \in E_{j,k}$, then
\begin{equation}
\label{kerfajap}
(A_a-\epsilon) \langle F,F \rangle  \leq 
\sum_j \sum_{k}|\int K_{a^j}(x_{j,k},y) F(y)d\mu(y)|^2 \mu(E_{j,k}) \leq (B_a + \epsilon)  \langle F,F \rangle ,
\end{equation}
which means
\begin{equation}
\label{mtvest}
(A_a-\epsilon) \langle F,F \rangle  \leq 
\sum_j \sum_{k} |(F,\phi_{j,k})|^2 \leq (B_a + \epsilon)  \langle F,F \rangle ,
\end{equation}
where
\begin{equation}
\label{phjkdf}
\phi_{j,k}(x)=[\mu(E_{j,k})]^{1/2} \overline{K_{a^j}}(x_{j,k},x).
\end{equation}

In our earlier article \cite{gm2}, we showed that (\ref{mtvest}) indeed holds,
provided the $E_{j,k}$ are also ``not too small'' (precisely, if they satisfy (\ref{measgeq}) directly below).
In fact, in Theorem \ref{framainfr} of that article, we showed (a more general form of) the following result:

\begin{theorem}
\label{framain}
Fix $a >1$, and say $c_0 , \delta_0 > 0$.  Suppose $f \in {\mathcal S}(\RR^+)$, 
and $f(0) = 0$.  Suppose that the Daubechies condition $(\ref{daub})$ holds.
Then there exists a constant $C_0 > 0$ $($depending only on ${\bf M}, f, a, c_0$ and $\delta_0$$)$
as follows:\\
For $t > 0$, let $K_t$ be the kernel of $f(t^2\Delta)$.
Say $0 < b < 1$.
Suppose that, for each $j \in \ZZ$, we can write ${\bf M}$ as a finite disjoint union of measurable sets 
$\{E_{j,k}: 1 \leq k \leq N_j\}$,
where:
\begin{equation}
\label{diamleq}
\mbox{the diameter of each } E_{j,k} \mbox{ is less than or equal to } ba^j, 
\end{equation}
and where:
\begin{equation}
\label{measgeq}
\mbox{for each } j \mbox{ with } ba^j < \delta_0,\: \mu(E_{j,k}) \geq c_0(ba^j)^n.
\end{equation}
$($In \cite{gm2} we show that such $E_{j,k}$ exist 
provided $c_0$ and $\delta_0$ are sufficiently small, independent of the 
values of $a$ and $b$.$)$

For $1 \leq k \leq N_j$, define $\phi_{j,k}$ by (\ref{phjkdf}).  Then
if $P$ denotes the projection in $L^2({\bf M})$ onto the space of constants, we have
\[(A_a-C_0b) \langle F,F \rangle  \leq 
\sum_j \sum_{k} |(F,\phi_{j,k})|^2 \leq (B_a + C_0b)  \langle F,F \rangle ,\]
for all $F \in (I-P)L^2({\bf M})$.  In particular, if $A_a - C_0b > 0$, then
$\left\{\phi_{j,k}\right\}$ is a frame for $(I-P)L^2({\bf M})$,
with frame bounds $A_a - C_0b$ and $B_a + C_0b$.
\end{theorem}  

Thus, in these circumstances, if $b$ is sufficiently small,
$\{\phi_{j,k}\}$ is a frame, in fact a smooth, nearly tight frame, since
\[ \frac{B_a + C_0b}{A_a - C_0b} \sim \frac{B_a}{A_a} = 
1 + O(|(a-1)^2 (\log|a-1|)|). \]

To justfiy the formal argument leading from (\ref{aasumba0}) to (\ref{phjkdf}),
and to go beyond the $L^2$ theory, one needs the following information about the kernel
$K_t$, which we established in Lemma \ref{manmolcw} of our earlier paper \cite{gmcw} (see also the remark following 
the proof of that lemma):

\begin{lemma}
\label{manmol}
Say $f(0) = 0$.  Then
for every pair of
$C^{\infty}$ differential operators $X$ $($in $x)$ and $Y$  $($in $y)$ on ${\bf M}$,
and for every integer $N \geq 0$, there exists $C_{N,X,Y}$ as follows.  Suppose
$\deg X = j$ and $\deg Y = k$.  Then
\begin{equation}
\label{diagest}
t^{n+j+k} \left|\left(\frac{d(x,y)}{t}\right)^N XYK_t(x,y)\right| \leq C_{N,X,Y} 
\end{equation}
for all $t > 0$ and all $x,y \in {\bf M}$.  (The result holds even without the hypothesis
that $f(0) = 0$, provided we look only at $t \in (0,1]$.)
\end{lemma}

The main results are Theorems \ref{besmain2} and \ref{besmain3} below, whose precise statements can be read now.
To summarize them: fix any $M_0 > 0$.
We study frame expansions for the space $B_{p,0}^{\alpha q}$, consisting
of distributions $F$ in the Besov space $B_{p}^{\alpha q}$ on ${\bf M}$ for which $F1 = 0$.
We assume that the $E_{jk}$ satisfy the conditions of Theorem \ref{framain} above
for $b$ sufficiently small, and also that, if $0 < p < 1$, and if
$ba^j \geq \delta_0$, then $\mu(E_{j,k}) \geq {\cal C}$ (for some 
${\cal C} > 0$).  (Such sets $E_{j,k}$ are easily constructed.)  We assume that $f(s) = s^l f_0(s)$ for some $f_0
\in {\cal S}(\RR^+)$, and for $l$ sufficiently large, depending on the indices $p,q, \alpha$
(so that $f(t^2\Delta) = t^{2l}\Delta^l f_0(t^2\Delta)$).
We let the 
$\phi_{j,k}$ be as in Theorem \ref{framain}, and let 
$\varphi_{j,k}(x)=  \phi_{j,k}(x)/[\mu(E_{j,k})]^{1/2} = \overline{K_{a^j}}(x_{j,k},x)$.   We then show that
a distribution $F$, of order at most $M_0$ and satisfying $F1 = 0$, is in $B_{p,0}^{\alpha q}$ if and only if 
\[ (\sum_{j = -\infty}^{\infty} a^{-j\alpha q} [\sum_k \mu(E_{j,k})
|\langle  F, \varphi_{j,k}  \rangle  |^p]^{q/p})^{1/q} < \infty;\]
further this
expression furnishes a norm on $B_{p,0}^{\alpha q}$ which is equivalent to the usual norm.  Moreover,
if $F \in B_{p,0}^{\alpha q}$, there exist constants $r_{j,k}$ with 
\begin{equation}
\label{rjknrmdf}
(\sum_{j = -\infty}^{\infty} a^{-j\alpha q} [\sum_k \mu(E_{j,k})|r_{j,k}|^p]^{q/p})^{1/q} < \infty 
\end{equation}
such that 
\begin{equation}
\label{besexpF}
F = \sum_{j=-\infty}^{\infty} \sum_k \mu(E_{j,k}) r_{j,k} \varphi_{j,k}
\end{equation}
with convergence in $B_p^{\alpha q}$;  and the infimum of the sums in (\ref{rjknrmdf}), taken over
all collections of numbers $\{r_{j,k}\}$ for which (\ref{besexpF}) holds, defines a norm on 
$B_{p,0}^{\alpha q}$, which is equivalent to the usual norm.\\

In addition to Lemma \ref{manmol},
our main tools will be the characterization of Besov spaces on $\RR^n$ by Frazier and Jawerth
\cite{FJ1},  and the characterization of these spaces on ${\bf M}$ by Seeger and 
Sogge \cite{SS}.  Our methods are largely adapted from those of Frazier and Jawerth.
There are, however, at least three major differences:\\
\ \\
1. We need to find replacements, on ${\bf M}$, for the condition
that a function have numerous vanishing moments.  Specifically, note that if $g \in C_c^{\infty}(\RR^n)$,
then $\Delta^l g$ has $2l-1$ vanishing moments for any $l \geq 1$, if $\Delta$ is the usual Laplacian.
In order to make effective use of our frames in Besov spaces,
we need an analogue of this on ${\bf M}$.  Say $g \in C^{\infty}({\bf M})$; what replacement condition
does $\Delta^l g$ satisfy, if, as usual, $\Delta$ is the Laplace-Beltrami operator on ${\bf M}$?
We will find an effective replacement in Lemma \ref{fjan2} of the next section.  These considerations 
explain why we need to 
use functions $f$ of the form $f(s) = s^l f_0(s)$ for $l$ sufficiently large, depending on the indices.\\  
\ \\
2. If one knows appropriate information
about the size of the frame coefficients of a function $F$, then, by adapting the methods of 
Frazier-Jawerth and by using the results of Seeger-Sogge, we learn only that 
$SF$ (not $F$) is in the desired Besov space, where
\begin{equation}
\label{sover}
SF = \sum_j \sum_k \mu(E_{j,k})  \langle  F,\phi_{j,k}  \rangle  \phi_{j,k};
\end{equation}
another step is then required.  Although, for $b$ sufficiently small, $SF$ is an excellent approximation 
to a multiple of $F$ (since the $\{\phi_{j,k}\}$ are a nearly tight frame),
it generally does not {\em equal} a multiple of $F$.  To conclude that $F$ itself is in the desired Besov space, 
we will need to use the theory of pseudodifferential operators (in Theorem \ref{besmain1} below).\\
\ \\
3. We need to show that $S$ is bounded on the Besov spaces, and a technical issue arises when the index $p$
lies between $0$ and $1$.  Since the $p$-triangle inequality generally becomes more and more wasteful when one splits 
quantities more finely (e.g. if we write $x > 0$ as $\frac{x}{N} + \ldots + \frac{x}{N}$ ($N$ terms), then
we find the wasteful estimate $x^p = (\frac{x}{N} + \ldots + \frac{x}{N})^p \leq N\frac{x^p}{N^p} = N^{1-p}x^p$),
and since we must use fine grids (that is, we must take $b$ to be small),
a rather subtle ``regrouping'' (or ``amalgamation'') argument is
needed at one point (Theorem \ref{sumopbes} below).\\
\ \\
In order for the notation to be fully analogous to that in \cite{FJ1}, we shall need
to adapt the notation that we used in our earlier article \cite{gm2};
through much of this article, we will write $E^j_k = E_{-j,k}$,
$x^j_k = x_{-j,k}$, $\varphi^j_k = \varphi_{-j,k}$. (The notations on the right sides of these
equations were used frequently in \cite{gm2}.)

\subsection{Historical Comments}

Although we are adapting the methods of Frazier and Jawerth \cite{FJ1}, we should note that
they were not working with nearly tight frames, but rather with the $\varphi$-transform.  Characterizations
of Besov spaces on $\RR^n$, which are similar to ours,
were obtained by Gr\"ochenig \cite{groch} (see also \cite{fg0}, \cite{fg1}, \cite{fg2}) through use of frames, and by
Meyer \cite{meyer}, through use of bases of orthonormal wavelets.

In \cite{dahsch}, Dahmen and Schneider 
used parametric liftings from standard bases on the unit cube
 to obtain biorthogonal wavelet bases on manifolds which are the disjoint union of smooth parametric images of 
the standard cube.  Using these bases, they obtained characterizations of the Besov spaces 
$B_p^{\alpha q}({\bf M})$, for $0 < p \leq \infty$, $q \geq 1$, and $\alpha > 0$.  Their results 
hold on manifolds with less than $C^{\infty}$ regularity (for a range of $\alpha$);
also, they applied their methods to the discretization of elliptic operator equations.  We consider
neither of these topics here.  However, our methods have the advantage of holding for all ${\bf M}$, 
and all $p,q, \alpha$.  Our frames have the advantage of being nearly tight, and admitting a space-frequency
analysis.  Moreover our results are coordinate-free, in the sense that our frames are constructed
without patching the manifold with charts.  We presume that this would lead to greater stability in applications,
if data is moving around the manifold in time, since one does not have to worry about data moving from
chart to chart, although this presumed advantage has not been established. 

In \cite{narc2}, Narcowich, Petrushev and Ward obtain a characterization of both Besov and Triebel-Lizorkin
spaces, through the size of frame coefficients, in the special case ${\bf M} = S^n$.  As frames they use
the ``needlets'' that they constructed in \cite{narc1}.  We discussed the similarities, advantages and
disadvantages of these frames as compared to ours on $S^n$, in section 3 of our earlier article \cite{gm2}.
They proved and used a result similar to our Lemma \ref{manmol}, and our methods (based on adapting the ideas in
\cite{FJ1}) are rather similar to theirs.  However, on the sphere, they constructed tight frames, so they did not need to 
deal with the issues \#1,2 and 3 above. 

Han and Sawyer \cite{hansaw} define Besov spaces on general spaces of homogeneous type, for a range of 
indices.  In \cite{hanyang}, Han and Yang give a characterization of these spaces using frames which they
construct.  These frames cannot be expected to be nearly tight, nor (on ${\bf M}$) have they been shown
to admit a space-frequency analysis.  Further, in the very general situation of \cite{hanyang}, there are no
derivatives, so results are obtained there only for smoothness index $\alpha \in (0,1)$.

\section{Integrating Products}

We shall need the following basic facts, from section 3 of \cite{gmcw}, about ${\bf M}$ and its
geodesic distance $d$.
For $x \in {\bf M}$, we let $B(x,r)$ denote the ball $\{y: d(x,y) < r\}$.
\begin{proposition}
\label{ujvj}
Cover ${\bf M}$ with a finite collection of open sets $U_i$  $(1 \leq i \leq I)$,
such that the following properties hold for each $i$:
\begin{itemize}
\item[$(i)$] there exists a chart $(V_i,\phi_i)$ with $\overline{U}_i
\subseteq V_i$; and
\item[$(ii)$] $\phi_i(U_i)$ is a ball in $\RR^n$.
\end{itemize}
Choose $\delta > 0$ so that $3\delta$ is a Lebesgue number for the covering
$\{U_i\}$.  Then, there exist $c_1, c_2 > 0$ as follows:\\
For any $x \in {\bf M}$, choose any $U_i \supseteq B(x,3\delta)$.  Then, in
the coordinate system on $U_i$ obtained from $\phi_i$, 
\begin{equation}
\label{rhoeuccmp2}
d(y,z) \leq c_2|y-z|
\end{equation}
for all $y,z \in U_i$; and 
\begin{equation}
\label{rhoeuccmp}
c_1|y-z| \leq d(y,z) 
\end{equation}
for all $y,z \in B(x,\delta)$.
\end{proposition}
We fix collections $\{U_i\}$, $\{V_i\}$, $\{\phi_i\}$ and also $\delta$ as in 
Proposition \ref{ujvj}, once and for all.
\begin{itemize}
\item Notation as in Proposition \ref{ujvj}, 
there exist $c_3, c_4 > 0$, such that, whenever $x \in {\bf M}$ and $0 < r \leq \delta$,
\begin{equation}
\label{ballsn}
c_3r^n \leq \mu(B(x,r)) \leq c_4r^n
\end{equation}
and such that, whenever $x \in {\bf M}$ and $r > \delta$,
\begin{equation}
\label{ballsn1}
c_3 \delta^n \leq \mu(B(x,r)) \leq \mu({\bf M}) \leq c_4r^n.
\end{equation}
\item 
For any $N > n$ there exists $C_N$ such that, for all $x \in {\bf M}$ and $t > 0$,
\begin{equation}
\label{ptestm}
\int_{\bf M} [1 + d(x,y)/t]^{-N} d\mu(y) \leq C_N t^n.
\end{equation}
\item 
  For any $N > n$ there exists $C_N^{\prime}$ such that, for all $x \in {\bf M}$ and $t > 0$,
\begin{equation}
\label{ptestm1}
\int_{d(x,y) \geq t} d(x,y)^{-N} d\mu(y) \leq C_N^{\prime} t^{n-N}.
\end{equation}
\end{itemize}

In Lemma 3.3 of \cite{FJ1}, Frazier and Jawerth proved, in essence, the following key
lemma on $\RR^n$, 
for which we must find analogues on ${\bf M}$.

\begin{lemma}
\label{fjrn}
Say $L, M$ are integers with $L \geq -1$ and $M \geq L+n+1$.  Then there exists $C > 0$ as follows.

Supppose $\varphi_1 \in C(\RR^n)$ and $\varphi_2 \in C^{L+1}(\RR^n)$ satisfy, for some $\sigma, \nu \in \ZZ$ with
$\sigma \geq \nu$, 
\[ |\varphi_1(x)| \leq (1+2^{\sigma}|x|)^{-M}, \]
\begin{equation}
\label{momfjrn}
\int x^{\alpha} \varphi_1(x) dx = 0 \mbox{  whenever } |\alpha| \leq L,
\end{equation}
and
\[ |\partial^{\gamma}\varphi_2(x)| \leq (1+2^{\nu}|x|)^{L+n+1-M} \mbox{  whenever } |\gamma| = L+1. \]
Then
\[ |(\varphi_1 * \varphi_2)(x)| \leq C 2^{-\sigma(L+n+1)}(1+2^{\nu}|x|)^{L+n+1-M}. \]
\end{lemma}
To clarify, (\ref{momfjrn}) is the empty condition if $L = -1$.

We will need two different analogues of this lemma.  The first is a straightforward adaptation of
Lemma (\ref{fjrn}) and its proof.

\begin{lemma}
\label{fjan1}
Say $L, M$ are integers with $L \geq -1$ and $M \geq L+n+1$.  Then there exists $C > 0$ as follows.

Say $\sigma \in \ZZ$, $\nu \in \RR$ with $2^{\sigma} \geq a^{\nu}$.

Say $x_0 \in {\bf M}$.  Select one of the charts $U_i$ (as in Proposition \ref{ujvj}) with
$B(x_0,3\delta) \subseteq U_i$.  Suppose, that in local coordinates on $U_i$, $Q$
is a dyadic cube of side $2^{-\sigma}$, and $3Q \subseteq B(x_0,\delta)$.   (Here $3Q$ is the cube
with the same center as $Q$ and $3$ times the side length $l(Q)$ of $Q$.)

Suppose also that 
 $\varphi_1 \in C({\bf M})$ satisfies the following conditions:
\[ \supp \varphi_1 \subseteq 3Q; \]
\[ |\varphi_1(y)| \leq 1 \mbox{  for all } y \in {\bf M}; \]
\begin{equation}
\label{momfjm}
\int y^{\alpha} \varphi_1(y) d\mu(y) = 0 \mbox{  whenever } |\alpha| \leq L,
\end{equation}
Also suppose $x_1 \in {\bf M}$, that $\varphi_2 \in C^{L+1}({\bf M})$, and that for all $y \in 3Q$,
\[ |\partial^{\gamma}\varphi_2(y)| \leq (1+a^{\nu}d(y,x_1))^{L+n+1-M} \mbox{  whenever } |\gamma| = L+1. \]
Then, 
\[ |\int_{\bf M}(\varphi_1 \varphi_2)(y)d\mu(y)| \leq C 2^{-\sigma(L+n+1)}(1+a^{\nu}d(x_0,x_1))^{L+n+1-M}. \]
\end{lemma}
To clarify, (\ref{momfjm}) is the empty condition if $L = -1$.\\
{\bf Proof} In this proof, $C$ will always denote a constant which depends only on $L$ and $M$ (and $n$, $a$ and
$({\bf M},g)$, of course); it may change from one line to the next.

Surely $|\int_{\bf M}(\varphi_1 \varphi_2)(y)d\mu(y)| = |\int_{3Q}(\varphi_1 \varphi_2)(y)d\mu(y)|$.
We write 
\begin{eqnarray*}
|\int_{3Q}(\varphi_1 \varphi_2)(y)d\mu(y)| & = &
|\int_{3Q}\varphi_1(y)[\varphi_2(y) - 
\sum_{|\gamma| \leq L} \frac{\partial^{\gamma} \varphi_2(x_0)(y-x_0)^{\gamma}}{\gamma!}] d\mu(y)|\\
& \leq & C(\int_{\{y \in 3Q: d(y,x_0) \leq c_1d(x_0,x_1)/2c_2\}} + 
\int_{\{y \in 3Q: d(y,x_0) > c_1d(x_0,x_1)/2c_2\}})|x_0-y|^{L+1} \Phi(y) d\mu(y)\\
& = & I + II,
\end{eqnarray*}
where $c_1$ and $c_2$ are as in Proposition \ref{ujvj}, and
where $\Phi(y) = \sup_{0 < \epsilon < 1} \sum_{|\gamma|=L+1}|\partial^{\gamma}\varphi_2(x_0 + \epsilon(y-x_0))|$.
In I, we have that whenever $0 < \epsilon < 1$, then
\begin{equation}
\label{linsegest}
d(x_0 + \epsilon(y-x_0),x_0) \leq c_2|(x_0 + \epsilon(y-x_0)-x_0| \leq c_2|y-x_0| \leq (c_2/c_1)d(y,x_0) \leq
d(x_0,x_1)/2,
\end{equation}
so that in I, by the hypotheses,
\[ |\Phi(y)| \leq C (1+a^{\nu}d(x_0,x_1))^{L+n+1-M}. \]
In II we just note that
\[ |\Phi(y)| \leq 1. \]
Accordingly
\[ I \leq C (1+a^{\nu}d(x_0,x_1))^{L+n+1-M}\int_{3Q} |x_0-y|^{L+1} d\mu(y) \leq 
C 2^{-\sigma(L+n+1)}(1+a^{\nu}d(x_0,x_1))^{L+n+1-M}, \]
while if $a^{\nu}d(x_0,x_1) \leq 1$, we have
\[ II \leq \int_{3Q} |x_0-y|^{L+1}d\mu(y) \leq C 2^{-\sigma(L+n+1)}
\leq C 2^{-\sigma(L+n+1)}(1+a^{\nu}d(x_0,x_1))^{L+n+1-M}. \]
Finally, say $a^{\nu}d(x_0,x_1) > 1$.  Note that $1 \leq C[2^{\sigma}|x_0-y|]^{-1}$ for all $y \in 3Q$.  Raising this
to the $M$th power, and then using (\ref{ptestm1}) and the assumption that $2^{\sigma} \geq a^{\nu}$, we see that
\begin{eqnarray} 
II & \leq & C  2^{-\sigma M}\int_{\{y \in 3Q: d(y,x_0) > c_1 d(x_0,x_1)/2c_2\}}|x_0-y|^{L-M+1}d\mu(y)\nonumber\\
& \leq & C  2^{-\sigma M}\int_{d(y,x_0) > c_1 d(x_0,x_1)/2c_2}d(y,x_0)^{L-M+1}d\mu(y)\label{IIstrt}\\
& \leq & C  2^{-\sigma M} d(x_0,x_1)^{L+n+1-M} \nonumber\\
& = & C 2^{-\sigma(L+n+1)} [2^{\sigma} d(x_0,x_1)]^{L+n+1-M} \nonumber\\
& \leq & C 2^{-\sigma(L+n+1)}(1+a^{\nu}d(x_0,x_1))^{L+n+1-M},\label{IIend} 
\end{eqnarray}
as claimed.\\

In our second analogue of Lemma \ref{fjrn}, instead of assuming that $\varphi_1$
is supported in a chart and satisfies familiar moment conditions there, as we did in Lemma \ref{fjan1},
we will instead allow $\varphi_1$ to be supported anywhere in ${\bf M}$.  The moment conditions will be
replaced by an assumption that $\varphi  = \Delta^l \Phi$ for another well-behaved function $\Phi$.
Formally, in this lemma,
the role of $L$ in Lemma \ref{fjan1} will be played by $2l-1$.
\begin{lemma}
\label{fjan2}
Say $l, M$ are integers with $l \geq 0$ and $M > n$.  
Then there exists $C > 0$ as follows.

Say $\sigma, \nu \in \RR$ with $\sigma \geq \nu$.

Say $x_0 \in {\bf M}$, and suppose that 
 $\varphi_1 = \Delta^l\Phi$, where $\Phi \in C^{2l}({\bf M})$ satisfies:
\[ |\Phi(y)| \leq (1+a^{\sigma}d(y,x_0))^{-M}. \]
Also suppose $x_1 \in {\bf M}$, that $\varphi_2 \in C^{2l}({\bf M})$, and that for all $y \in {\bf M}$,
\[ |\Delta^l \varphi_2(y)| \leq (1+a^{\nu}d(y,x_1))^{n-M}.\] 
Then, 
\[ |\int_{\bf M}(\varphi_1 \varphi_2)(y)d\mu(y)| \leq C a^{-\sigma n}(1+a^{\nu}d(x_0,x_1))^{n-M}. \]
\end{lemma}
{\bf Proof} In this proof, $C$ will always denote a constant which depends only on $l, M$ and $n$ (and 
$a, ({\bf M},g)$, of course); it may change from one line to the next.

First we observe that we may take $l = 0$.  Indeed, if this case were known, then we could
prove the general result simply by noting that $\Delta$ is self-adjoint, so that

\[ \int_{\bf M}(\varphi_1 \varphi_2)(y)d\mu(y) = \int_{\bf M}(\Delta^l \Phi \varphi_2)(y)d\mu(y) =
\int_{\bf M}(\Phi \Delta^l \varphi_2)(y)d\mu(y). \]

So we may assume $l=0$.

We have
\[ |\int_{\bf M}(\varphi_1 \varphi_2)(y)d\mu(y)|  \leq 
 I+II, \]
where we are setting
\[ I = \int_{d(y,x_0) \leq d(x_0,x_1)/2}|\varphi_1(y)\varphi_2(y)| d\mu(y), \]
\[ II = \int_{d(y,x_0) > d(x_0,x_1)/2}|\varphi_1(y)\varphi_2(y)| d\mu(y). \]

We shall show that I and II are less than or equal to $C a^{-\sigma n}(1+a^{\nu}d(x_0,x_1))^{n-M}$,
and then we will be done.  

For I and II we need to estimate $\varphi_2(y)$.   We shall use the evident estimates:
\begin{equation}
\label{estIphi}
\mbox{In I, } |\varphi_2(y)| \leq (1+a^{\nu}d(x_0,x_1))^{n-M},
\end{equation}
and 
\begin{equation}
\label{estIIphi}
\mbox{In II, } |\varphi_2(y)| \leq 1.
\end{equation}

From (\ref{estIphi}) and the hypotheses on $\varphi_1 = \Phi$, we find that
\[ I \leq 
(1+a^{\nu}d(x_0,x_1))^{n-M}\int_{\bf M} (1+a^{\sigma}d(y,x_0))^{-M}d\mu(y) \leq
 Ca^{-\sigma n}(1+a^{\nu}d(x_0,x_1))^{n-M} \]
as needed.  (We have used (\ref{ptestm}).)

For II, we have the estimate
\[ II \leq \int_{d(y,x_0) > d(x_0,x_1)/2}(1+a^{\sigma}d(y,x_0))^{-M}d\mu(y)\]
Suppose first that $a^{\nu}d(x_0,x_1) \leq 1$.  Then we can just note that, by (\ref{ptestm}),
\[ II \leq C a^{-\sigma n } \leq C a^{-\sigma n}(1+a^{\nu}d(x_0,x_1))^{n-M}. \]

If, instead, $a^{\nu}d(x_0,x_1) > 1$, we find that
\begin{eqnarray*} 
II & \leq & C  a^{-\sigma M}\int_{d(y,x_0) > d(x_0,x_1)/2}d(y,x_0)^{-M}d\mu(y)\\
& \leq & C a^{-\sigma n} (1+a^{\nu}d(x_0,x_1))^{n-M} 
\end{eqnarray*}
as an argument just like the one
beginning with (\ref{IIstrt}) and ending with (\ref{IIend}) shows.  This completes the proof.\\

Next we need analogues of Lemma 3.4 of \cite{FJ1}.  After one multiplies by certain constants,
that lemma states:

\begin{lemma}
\label{fjrn2}
If $Q$ is a dyadic cube of $\RR^n$, let $x_Q$ denote its center.

Let $1 \leq p \leq \infty$, and suppose $\sigma, \eta \in \ZZ$, $\eta \leq \sigma$.  Suppose $F(x) = 
\sum_{l(Q) = 2^{-\sigma}} s_Q f_Q(x)$, where 
\[ |f_Q(x)| \leq 2^{-\sigma n}(1+2^{\eta}|x-x_Q|)^{-n-1}. \]
Then
\[ \|F\|_{L^p} \leq C 2^{-\eta n} (\sum_{l(Q) = 2^{-\sigma}} 2^{-\sigma n}|s_Q|^p)^{1/p}. \]
\end{lemma}

Here is our first analogue.

\begin{lemma}
\label{fjan3}
Let $1 \leq p \leq \infty$, and as usual fix $a > 1$.  Also fix $b > 0$.
Then there exists $C > 0$ as follows.

Suppose $\eta \in \RR$, $j \in \ZZ$, $\eta \leq j$. 
Write ${\bf M}$ as a finite disjoint union of 
measurable subsets $\{E^j_k: 1 \leq k \leq {\cal N}_j\}$, each of diameter less than $ba^{-j}$.
For each $k$ with $1 \leq k \leq {\cal N}_j$, select any $x^j_k \in E^j_k$. 

Suppose that, for $x \in {\bf M}$,  $F(x) = 
\sum_{k=1}^{{\cal N}_j} s_{j,k} f_{j,k}(x)$, where 
\[ |f_{j,k}(x)| \leq \mu(E^j_k)(1+a^{\eta}d(x,x^j_k))^{-n-1}. \]
Then
\[ \|F\|_{L^p} \leq C a^{-\eta n} (\sum_{k} \mu(E^j_k)|s_{j,k}|^p)^{1/p}. \]
\end{lemma}
{\bf Proof} 
We have
\begin{equation}
\label{sjkamu}
\|F\|_p \leq (\int_{\bf M} 
[\sum_{k=1}^{{\cal N}_j} |s_{j,k}| \mu(E^j_k)(1+a^{\eta}d(x,x^j_k))^{-n-1}]^p d\mu(x))^{1/p}
\end{equation}
Let $Y_j$ be the finite measure space $\{1,\ldots,{\cal N}_j\}$ with measure $\lambda$
where $\lambda(\{k\}) = \mu(E^j_k)$ for each $k \in Y_j$. and define 
${\cal K}: L^p(Y_j) \rightarrow L^p({\bf M})$ by ${\cal K} r(x) = \int_{Y_j} K(x,k)r(k) d\lambda(k)$ for $r \in L^p(Y_j)$,
where
\[ K(x,k) = (1+a^{\eta}d(x,x^j_k))^{-n-1}. \]
By standard arguments, $\|{\cal K} r\|_p \leq M\|r\|_p$, where
$M$ is any number satisfying
\begin{equation}
\label{mkxkmu}
M \geq \int_{\bf M} K(x,k) d\mu(x)
\end{equation}
for all $k \in Y_j$ and also
\begin{equation}
\label{mkxkla}
M \geq \int_{Y_j} K(x,k) d\lambda(k)
\end{equation}
for all $x \in {\bf M}$.
By (\ref{sjkamu}), $\|F\|_p \leq \|{\cal K}r\|_p$ if $r(k) = |s_{j,k}|$.
Thus we need only show that we may take $M = C a^{-\eta n}$.  But, for this $M$, (\ref{mkxkmu}) 
holds by (\ref{ptestm}).  As for (\ref{mkxkla}), choose $B > \max(2b,1)$.  If $x \in {\bf M}$, we have
\begin{eqnarray}
\int_{Y_j} |K(x,k)| d\lambda(k) & \leq & C\sum_{k} \mu(E^j_k)(B+a^{\eta}d(x,x^j_k))^{-n-1} \label{m2stt}\\ 
& \leq & C \sum_{k} \int_{E^j_k}(B+a^{\eta}d(x,y))^{-n-1} d\mu(y) \nonumber \\
& \leq & C\int_{\bf M}(1+a^{\eta}d(x,y))^{-n-1} d\mu(y) \nonumber \\
& \leq & Ca^{-\eta n} \label{m2end}
\end{eqnarray}
as desired.  Here we have used the assumption that $\eta \leq j$, the fact that the diameter
of $E^j_k$ is at most $ba^{-j}$, and (\ref{ptestm}).  This completes the proof.\\

Our second analogue deals only with $L^p$ norms of functions defined on finite sets.
\begin{lemma}
\label{fjan4}
Let $1 \leq p \leq \infty$, and as usual fix $a > 1$.  Also fix $b > 0$.
Then there exists $C > 0$ as follows.

Say $j \in \ZZ$.  Select sets $E^j_k$ and points $x^j_k$ as in Lemma \ref{fjan3}.  Let $Y_j$ 
again be the finite measure space $\{1,\ldots,{\cal N}_j\}$ with measure $\lambda$
where $\lambda(\{k\}) = \mu(E^j_k)$ for each $k \in Y_j$. 

Suppose $U_i$ is one of the charts of Proposition \ref{ujvj}.  Say $B(x_0,3\delta) \subseteq U_i$.
Say $\sigma \in \ZZ$, and let
${\cal Q}_{\sigma}$ denote the set of dyadic cubes of length 
$2^{-\sigma}$  in $U_i$ which are contained in $B(x_0,\delta)$.  (We are using local coordinates on 
$U_i$.) \\
(a) Say $2^{\sigma} \geq a^{j}$.  Suppose that, for $k \in Y_j$,  $F(k) = 
\sum_{Q \in {\cal Q}_{\sigma}} s_Q f_{Q}(k)$, where 
\[ |f_Q(k)| \leq 2^{-\sigma n}(1+a^{j}d(x_Q,x^j_k))^{-n-1}. \]
Then
\[ \|F\|_{L^p(Y_j)} \leq C a^{-jn}(\sum_{l(Q) = 2^{-\sigma}} 2^{-\sigma n}|s_{Q}|^p)^{1/p}. \]
(b) Suppose instead $2^{\sigma} \leq a^{j}$.  Suppose that, for $k \in Y_j$,  $F(k) = 
\sum_{Q \in {\cal Q}_{\sigma}}  s_Q f_{Q}(k)$, where
\[ |f_Q(k)| \leq 2^{-\sigma n}(1+2^{\sigma}d(x_Q,x^j_k))^{-n-1}. \]
Then
\[ \|F\|_{L^p(Y_j)} \leq C 2^{-\sigma n}(\sum_{l(Q) = 2^{-\sigma}} 2^{-\sigma n}|s_{Q}|^p)^{1/p}. \]
\end{lemma}
{\bf Proof} Say ${\cal Q}_{\sigma} = \{Q_1,\ldots,Q_{I(\sigma)}\}$, and let
$X_{\sigma} = \{1,\ldots,I(\sigma)\}$, with measure $\tau$, where $\tau(m) = 2^{-\sigma n}$
for each $m \in X_{\sigma}$. 

In either (a) or (b), we define an operator
${\cal K}: L^p(X_{\sigma}) \rightarrow L^p(Y_j)$ by
${\cal K} r(k) = \int_{X_{\sigma}} K(k,m)r(m) d\tau(m)$ for $r \in L^p(X_{\sigma})$,
where in (a),
\begin{equation}
\label{kkma}
K(k,m) = (1+a^{j}d(x_{Q_m},x^j_k))^{-n-1},
\end{equation}
while in (b)
\begin{equation}
\label{kkmb}
K(k,m) = (1+2^{\sigma}(x_{Q_m},x^j_k))^{-n-1}.
\end{equation}
With these definitions of $K(k,m)$, in either (a) or (b),
$\|F\|_p \leq \|{\cal K}r\|_p$ if $r(m) = |s_{Q_m}|$.
On the other hand, $\|{\cal K} r\|_p \leq M\|r\|_p$, 
where $M$ is any number satisfying
\begin{equation}
\label{mkkmla}
M \geq \int_{Y_j} K(k,m) d\lambda(k) 
\end{equation}
for all $m \in X_{\sigma}$ and also
\begin{equation}
\label{mkkmtau}
M \geq \int_{X_{\sigma}} K(k,m) d\tau(m)
\end{equation}
for all $k \in Y_j$.

In (a), where $K(k,m)$ is given by (\ref{kkma}), we need only show that we may take 
$M = C a^{-jn}$.  The fact that (\ref{mkkmla}) then holds follows just as in the argument starting with
(\ref{m2stt}) and ending with (\ref{m2end}).  
Similarly, for (\ref{mkkmtau}), say $k \in Y_j$.  Since $Q_m \subseteq B(x_0,\delta)$, the diameter of 
$Q_m$ is at most $c2^{-\sigma}$ for some $c$.  Choose $B > \max(2c,1)$.  Then
\begin{eqnarray*}
\int_{X_{\sigma}} K(k,m) d\tau(m) & \leq & C\sum_{m} 2^{-\sigma n}(B+a^{j}d(x_{Q_m},x^j_k))^{-n-1}\\ 
& \leq & C \sum_{m} \int_{Q_m}(B+a^{j}d(x,x^j_k))^{-n-1} d\mu(y)\\
& \leq & C B \int_{\bf M}(1+a^{j}d(x,x^j_k))^{-n-1} d\mu(y) \\
& \leq & C a^{-jn}
\end{eqnarray*}
as desired.  Here we have used the assumption that $2^{-\sigma} \leq a^{-j}$, the fact that the diameter
of $Q_m$ is at most $c2^{-\sigma}$, and (\ref{ptestm}).  

In (b), where $K(k,m)$ is given by (\ref{kkmb}), 
we need only show that we may take $M = C 2^{-\sigma n}$.
This, however, follows in a similar manner to (a), if one now uses the assumption that
$a^{-j} \leq 2^{-\sigma}$.   This completes the proof.

\section{Besov spaces}

We will need the following simple fact about operators on $l^q(\NN)$, for $0 < q \leq 1$,
which is again adapted from arguments in \cite{FJ1}.

\begin{proposition}
\label{lqconv}
Suppose $0 < q \leq \infty$.
Say $K: \NT \times \NT \rightarrow \RR$ is nonnegative.  If $z$ is a nonnegative
sequence, define the nonnegative sequence ${\cal K}z$ by
\[ ({\cal K}z)(r) = \sum_{s=0}^{\infty} K(r,s) z(s). \]
Let $M_q$ be a number satisfying
\[ M_q \geq [\sum_{s=0}^{\infty} K(r,s)^q]^{1/q} \]
for all $r$, and also
\[ M_q \geq [\sum_{r=0}^{\infty} K(r,s)^q]^{1/q} \]
for all $s$.  
Then:\\
(a) If $1 \leq q \leq \infty$, then 
 for every nonnegative sequence $z$, $\|{\cal K}z\|_q \leq M_1\|z\|_q$.\\
(b) If $0 < q < 1$, then 
for every nonnegative sequence $z$, $\|{\cal K}z\|_q \leq M_q\|z\|_q$.\\
(Here, $\|z\|_q$ denotes the $l^q(\NT)$ ``norm'' of $z$, which could be $\infty$;
and all nonnegative numbers here are allowed to be $\infty$.  Also, here and elsewhere,
we follow the usual rules for interpreting the expressions when $q=\infty$.)
\end{proposition}
{\bf Proof} (a) is of course well known.  For (ii), 
note that, by the $q$-triangle inequality,
\[ ({\cal K}z)(r)^q \leq \sum_{s=0}^{\infty} K(r,s)^q z(s)^q. \]
By the known case $q=1$ of the proposition, we now see that
\[ \|({\cal K}z)^q\|_1 \leq M_q^q \|z^q\|_1. \]
Raising both sides to the $1/q$ power, we obtain the desired result.\\

For the rest of this section, we fix $a > 1$.  We also fix $\alpha, p, q$ with
$-\infty < \alpha < \infty$ and $0 < p,q \leq \infty$. \\

We use the notation for inhomogeneous Besov spaces $B_p^{\alpha q}$ on $\RR^n$ from \cite{FJ1}.  
Thus, on $\RR^n$, one takes
any $\Phi \in {\cal S}$ supported in the closed unit ball, which does not vanish anywhere in the ball of
radius $5/6$ centered at $0$.  One also takes functions $\varphi_{\nu} \in {\cal S}$ for $\nu \geq 1$, 
supported in the annulus $\{\xi: 2^{\nu-1} \leq |\xi| \leq 2^{\nu+1}\}$, satisfying 
$|\varphi_{\nu}(\xi)| \geq c > 0$ for $3/5 \leq 2^{-\nu}|\xi| \leq 5/3$ and also 
$|\partial^{\gamma} \varphi_{\nu}| \leq c_{\gamma} 2^{-\nu \gamma}$ for every multiindex $\gamma$.
The Besov space $B_p^{\alpha q}(\RR^n)$ is then the space of $F \in {\cal S}'(\RR^n)$ such that

\[ \|F\|_{B_p^{\alpha q}} = \|\check{\Phi}*F\|_{L^p} + 
\left(\sum_{\nu = 0}^{\infty} (2^{\nu \alpha}\|\check{\varphi}_{\nu}*F\|_{L^p})^q \right)^{1/q} < \infty. \]

(Here we use the usual conventions if $p$ or $q$ is $\infty$.
The definition of $B_p^{\alpha q}(\RR^n)$ is independent of the choices of 
$\Phi, \varphi_{\nu}$ (\cite{Peet}, page 49).  Moreover, $B_p^{\alpha q}(\RR^n)$ is a quasi-Banach
space, and the inclusion $B_p^{\alpha q} \subseteq {\cal S}'$ is continuous
(\cite{Trieb}, page 48).  In particular the space $B^{\alpha}_{\infty,\infty}(\RR^n) =  {\cal C}^{\alpha}(\RR^n)$,
which is the usual H\"older space if $0 < \alpha < 1$, or in general a H\"older-Zygmund space for $\alpha > 0$ 
(\cite{Trieb}, page 51).  It is not hard to see, by using the definition and the Fourier transform, that
if $K \subseteq \RR^n$ is compact, and if $N$ is sufficiently large, then
\begin{equation}
\label{ccNbes}
\{F \in C^N: \mbox{ supp}F \subseteq K\} \subseteq B_p^{\alpha q}
\end{equation}
where the inclusion map is continuous if we regard the left side as a subspace of $C^N$.

Pseudodifferential operators of order $0$ are bounded on the Besov spaces (\cite{P}); in particular,
if $\psi \in C_c^{\infty}(\RR^n)$, the mapping $F \rightarrow \psi F$ is a bounded map on 
$B_p^{\alpha q}(\RR^n)$.   Moreover, if $\eta: \RR^n \rightarrow \RR^n$ is a diffeomorphism which equals
the identity outside a compact set, then one can define $F \circ \eta$ for $F \in B_p^{\alpha q}(\RR^n)$,
and the map $F \rightarrow F \circ \eta$ is bounded on the Besov spaces (\cite{Trieb}, chapter 2.10).
These facts then enable one
to define $B_p^{\alpha q}({\bf M})$: let $(W_i, \chi_i)$ be a finite atlas on ${\bf M}$ with charts $\chi_i$ 
mapping $W_i$ into the unit ball on $\RR^n$, and suppose $\{\zeta_i\}$ is a partition of unity subordinate
to the $W_i$.  Then one defines $B_p^{\alpha q}({\bf M})$ to be the space of distributions $F$ on 
${\bf M}$ for which 

\[ \|F\|_{B_p^{\alpha q}({\bf M})} = \sum_i \|(\zeta_i F) \circ \chi_i^{-1}\|_{B_p^{\alpha q}({\bf R}^n)} < \infty. \]

This definition does not depend on the choice of charts or partition of unity (\cite{Triebman}).\\

It will be convenient to fix a spanning set of the differential operators on ${\bf M}$ of degree less
than or equal to $J$ (for any fixed $J$).
Recall that we have already fixed a finite set ${\mathcal P}$ of real $C^{\infty}$ vector fields on ${\bf M}$, 
whose elements span the tangent space at each point.  For any integer $L \geq 1$, we let

\begin{equation}
\label{pmdfopdf}
{\mathcal P}^J = \{X_1\ldots X_M: X_1,\ldots,X_M \in {\mathcal P}, 1 \leq M \leq J\} \cup
\{\mbox{the identity map}\}.
\end{equation}
(In particular, ${\mathcal P}^1$ is what we have previously called ${\mathcal P}_0$.) \\

\begin{lemma}
\label{besov1}
Fix $b > 0$.  Also fix
an integer $l \geq 1$ with
\begin{equation}
\label{lgqgint}
2l > \max(n(1/p-1)_+ - \alpha, \alpha).
\end{equation}
where here $x_+ = \max(x,0)$.  Fix $M$ with $(M-2l-n)p > n+1$ if $0 < p < 1$, $M-2l-n > n+1$
otherwise.

Then there exists $C > 0$ as follows.

Say $j \in \ZZ$.  Select sets $E^j_k$ and points $x^j_k$ as in Lemma \ref{fjan3}.
Suppose that, for each $j \geq 0$, and each $k$, 
\begin{equation}
\label{var2lph}
\varphi^j_k = (a^{-2j}\Delta)^l\Phi^j_k,
\end{equation}
 where $\Phi^j_k \in C^{\infty}({\bf M})$ satisfies the following conditions:
\begin{equation}
\label{xphip}
|X\Phi^j_k(y)| \leq a^{j (\deg X + n)}(1+a^{j}d(y,x^j_k))^{-M} \mbox{  whenever } X \in {\mathcal P}^{4l}. 
\end{equation}
Then, for every $F$ in the inhomogeneous Besov space $B_p^{\alpha q}({\bf M})$, if we let 
\[ s_{j,k} = \langle F, \varphi^j_k \rangle, \]
then
\begin{equation}
\label{besov1way}
(\sum_{j = 0}^{\infty} a^{j\alpha q} [\sum_k \mu(E^j_k)|s_{j,k}|^p]^{q/p})^{1/q} \leq C\|F\|_{B_p^{\alpha q}}.
\end{equation}
\end{lemma}
{\bf Proof}
Cover ${\bf M}$ by a finite collection of open sets $\{W_r\}$, where each $W_r$ has the form
$B(x_r,\delta)$ for some $x_r \in {\bf M}$.  Let $\{\zeta_r\}$ be a partition of unity subordinate to 
the $\{W_r\}$.  Let $Z_r = \supp \zeta_r \subseteq W_r$.
Then \cite{Trieb} the map $F \rightarrow \zeta_r F$ is continuous from $B_p^{\alpha q}$ to 
itself.  Without loss, we may therefore assume that $\supp F \subseteq Z$ where $Z$ is a compact subset
of $W=B(x_0,\delta)$ for some $x_0$.  Choose a chart $U_i$, as in Proposition \ref{ujvj}, with 
$B(x_0,3\delta) \subseteq U_i$.  Select $\zeta \in C_c^{\infty}(W)$ with $\zeta \equiv 1$ 
in a neighborhood of $Z$.

In local coordinates on $U_i$, $U_i$ is a ball in $\RR^n$.  By \cite{FJ1} 
(changing their notation slightly by multiplying by certain constants) we may write 
\[ F = \sum_{m \in \ZZ^n} s_m b_m + \sum_{\sigma=0}^{\infty} \sum_{l(Q)=2^{-\sigma}} s_Q a_Q \]
with convergence in $B_p^{\alpha q}$.  Here, if $Q_{0m}$ is the dyadic cube of side $1$ with ''lower left corner'' $m$, 
\[ \mbox{supp} b_m \subseteq 3Q_{0m}, \]
\[ |\partial ^{\gamma} b_m| \leq 1 \mbox{  if } |\gamma| \leq 2l, \]
\[ \mbox{supp} a_Q \subseteq 3Q, \]
\[ |\partial ^{\gamma} a_Q| \leq |Q|^{-|\gamma|/n} \mbox{  if } |\gamma| \leq 2l, \]
\[ \int x^{\gamma} a_Q(x)dx = 0 \mbox{   if } |\gamma| \leq 2l-1,\]
and finally
\[ (\sum_{m \in \ZZ^n} |s_m|^p)^{1/p} +
 [\sum_{\sigma=0}^{\infty} 2^{\sigma \alpha q}[\sum_{l(Q)=2^{-\sigma}} 2^{-\sigma n}|s_Q|^p]^{q/p}]^{1/q} \leq C\|F\|_{B_p^{\alpha q}}. \]
Since $F = \zeta F$, we have
\begin{equation}
\label{fsmsq}
F =  \sum_{m \in \ZZ^n} s_m \zeta b_m + \sum_{\sigma=0}^{\infty} \sum_{l(Q)=2^{-\sigma}} s_Q \zeta a_Q 
\end{equation}
with convergence in $B_p^{\alpha q}$, hence in ${\cal E}'({\bf M})$.       
Now, the (Euclidean) distance from supp$\zeta$ to $W^c$ is positive. 
Thus there exists $\sigma_0$ with the property
that if $\sigma \geq \sigma_0$, $l(Q) = 2^{-\sigma}$ and $3Q \cap \mbox{supp}\zeta \neq \oslash$, then 
supp$\zeta + 3Q \subseteq W$.  (Note that $\sigma_0$ does not depend on $F$.)
Accordingly, if $l(Q) \leq 2^{-\sigma_0}$, then either $\zeta a_Q \equiv 0$ or $3Q \subseteq W$.  
Moreover, only finitely many cubes $3Q$ with $2^{-\sigma_0} < l(Q) \leq 1$ intersect the compact set supp$\zeta$;
let ${\cal Q}_0$ denote the collection of such cubes.  Thus we may write
\begin{equation}
\label{ff0sq}
\zeta F = F_0 + \sum_{\sigma=\sigma_0}^{\infty} \sum_{l(Q)=2^{-\sigma},3Q \subseteq W} s_Q \zeta a_Q
\end{equation}
where
\begin{equation}
\label{foqodf}
F_0 =  \sum_{m \in \ZZ^n, Q_{0m} \in {\cal Q}_0} s_m \zeta b_m +
\sum_{\sigma=0}^{\sigma_0-1} \sum_{l(Q)=2^{-\sigma}, Q \in {\cal Q}_0} s_Q \zeta a_Q 
\end{equation}

Let $c = \log_2 a$.  Then, since the series in (\ref{fsmsq}) converges to $F$ in ${\cal E}'({\bf M})$, 
\[ s_{j,k} = \langle F_0,\varphi^j_k \rangle +  
\sum_{\sigma_0 \leq \sigma \leq jc}\: \sum_{l(Q)=2^{-\sigma},3Q \subseteq W} s_Q \langle \zeta a_Q, \varphi^j_k \rangle
+ \sum_{\sigma > jc}\: \sum_{l(Q)=2^{-\sigma},3Q \subseteq W} s_Q \langle \zeta a_Q, \varphi^j_k \rangle. \]

For each $j \geq 0$, 
let $Y_j$ be the finite measure space $\{1,\ldots,{\cal N}_j\}$ with measure $\lambda$
where $\lambda(\{k\}) = \mu(E^j_k)$ for each $k \in Y_j$.  
For each $j \geq 0$, define $s_j: Y_j \rightarrow \CC$ by $s_j(k) = s_{j,k}$.  Also, for each $Q$ with 
$3Q \subseteq W$, define $u^Q_j: Y_j \rightarrow \CC$ by
 $u^Q_j(k) = \langle \zeta a_Q, \varphi^j_k \rangle$.  Finally, define $u^0_j: Y_j \rightarrow \CC$ by
$u^0_j(k) = \langle F_0,\varphi^j_k \rangle$.  Then

\begin{equation}
\label{sjsumal}
s_j = u^0_j +  
\sum_{\sigma_0 \leq \sigma \leq jc}\: \sum_{l(Q)=2^{-\sigma},3Q \subseteq W} s_Q u^Q_j
+ \sum_{\sigma > jc}\: \sum_{l(Q)=2^{-\sigma},3Q \subseteq W} s_Q u^Q_j.
\end{equation}

Define $h$ on $U_i$  by $d\mu = h dx$ there.  Note that if $3Q \subseteq W$ then $\int x^{\gamma} (a_Q/h) d\mu = 0$
if $|\gamma| \leq 2l-1$.

Now say that $\sigma > jc$, i.e., that $2^{\sigma} \geq a^{j}$.  Then by Lemma \ref{fjan1} 
(replacing ``$\nu$'' in that lemma by $j$, and with
$a_Q/h = \varphi_1$ and $h\zeta\varphi^j_k = a^{j(2l+n)}\varphi_2$), then
\begin{equation}
\label{uqjest1}
|u^Q_j(k)| = |\langle \zeta a_Q, \varphi^j_k \rangle| = |\langle  a_Q/h , h\zeta\varphi^j_k \rangle| 
\leq Ca^{j(2l+n)} 2^{-\sigma(2l+n)}(1+a^{j}d(x_Q,x^j_k))^{2l+n-M}.
\end{equation}

Say instead $\sigma_0 \leq \sigma \leq jc$, so that $2^{\nu} < a^{j}$.  Then by Lemma \ref{fjan2}
(replacing ``$\nu$'' in that lemma by 
$\sigma/c$, replacing ``$\sigma$'' in that lemma by $j$, 
and with $\Phi^j_k = a^{jn}\Phi$, $\varphi^j_k = a^{jn-2jl}\varphi_1$, and $\zeta a_Q = 2^{2\sigma l}\varphi_2$)
we have
\begin{equation}
\label{uqjest2}
|u^Q_j(k)| = |\langle \zeta a_Q, \varphi^j_k \rangle| \leq C2^{2\sigma l} a^{-2jl}(1+2^{\sigma}d(x_Q,x^j_k))^{n-M}. 
\end{equation}

Moreover, if $Q$ is one of the finitely many cubes in ${\cal Q}_0$, then again by Lemma \ref{fjan2}
(taking ``$\nu$'' in that lemma to be $0$, replacing ``$\sigma$'' in that lemma by $j$,
and with $\Phi^j_k = a^{jnl}\Phi$, $\varphi^j_k = a^{jn-2jl}\varphi_1$, 
and $C\zeta b_m = \varphi_2$ if $Q = Q_{0m}$ has side length $1$, or $C\zeta a_Q = \varphi_2$ otherwise) we have
\begin{equation}
\label{uqjest3}
|u^0_j(k)| \leq C{\cal T} a^{-2jl}(1+d(x_Q,x^j_k))^{n-M}, 
\end{equation}
where 
\[ {\cal T} = 
\sum_{m \in \ZZ^n, Q_{0m} \in {\cal Q}_0} |s_m| +
\sum_{\sigma=0}^{\sigma_0-1} \sum_{l(Q)=2^{-\sigma}, Q \in {\cal Q}_0} |s_Q|. \]

Say now $p \geq 1$, and let $\|\:\|_p$ denote $L^p(Y_j)$ norm.  From (\ref{sjsumal}) we obtain
\[ \|s_j\|_p \leq \|u^0_j\|_p +  
\sum_{\sigma_0 \leq \sigma \leq jc}\: \|\sum_{l(Q)=2^{-\sigma},3Q \subseteq W} s_Q u^Q_j\|_p
+ \sum_{\sigma > jc}\: \|\sum_{l(Q)=2^{-\sigma},3Q \subseteq W} s_Q u^Q_j\|_p.  \]

Let 
\[ A_j =  \|s_j\|_p =  [\sum_{k=1}^{{\cal N}_j} \mu(E^j_k)|s_{j,k}|^p]^{1/p},\:\:
B_{\sigma} = [\sum_{l(Q)=2^{-\sigma},3Q \subseteq W} 2^{-\sigma n}|s_Q|^p]^{1/p}. \]

Then by (\ref{uqjest1}), (\ref{uqjest2}), (\ref{uqjest3}) and Lemma \ref{fjan4}, we see that

\begin{equation}
\label{ajbnuest}
A_j \leq C({\cal T}a^{-2jl} + \sum_{\sigma_0 \leq \sigma \leq jc} a^{-2jl} 2^{2\sigma l}B_{\sigma}
+ \sum_{\sigma > jc} a^{2jl} 2^{-2\sigma l}B_{\sigma}).
\end{equation}

(In using Lemma \ref{fjan4}, we have noted that $3Q \subseteq W \Rightarrow Q \subseteq W$.)
Now also write

\[ A^{\alpha}_j = a^{j\alpha} A_j,\:\: B^{\alpha}_{\sigma} = 2^{\sigma \alpha} B_{\sigma}. \]

Then, by (\ref{ajbnuest}),

\begin{equation}
\label{ajbnalest}
A_j^{\alpha} \leq C({\cal T} a^{-j(2l-\alpha)}+ \sum_{\sigma=\sigma_0}^{\infty} K(j,\sigma) B^{\alpha}_{\sigma}),
\end{equation}

where 

\[ K(j,\sigma) = a^{-j(2l-\alpha)}2^{\sigma(2l-\alpha)} \mbox{    if } \sigma_0 \leq \sigma \leq jc, \]

\[ K(j,\sigma) = a^{j(2l+\alpha)}2^{-\sigma(2l+\alpha)} \mbox{    if } \sigma_0 \geq jc. \]

By (\ref{lgqgint}), $2l$ is more than $\max(\alpha,-\alpha)$.  
Recall also that $c = \log_2 a$.  Thus, by Proposition \ref{lqconv},
$\|(A_j^{\alpha})\|_q \leq C({\cal T}+\|(B_{\sigma}^{\alpha})\|_q)$.  Consequently
\begin{equation}
\label{bescmp1}
(\sum_{j = 0}^{\infty} a^{j\alpha q} [\sum_k \mu(E^j_k)|s_{j,k}|^p]^{q/p})^{1/q} 
\leq C[(\sum_{m \in \ZZ^n} |s_m|^p)^{1/p} +
(\sum_{\sigma=0}^{\infty} 2^{\sigma \alpha q}[\sum_{l(Q)=2^{-\sigma}} 2^{-\sigma n}|s_Q|^p]^{q/p})^{1/q}]
 \leq C\|F\|_{B_p^{\alpha q}}
\end{equation}
as desired, at least if $p \geq 1$.

If instead $0 < p < 1$, we evaluate each side of (\ref{sjsumal}) at $k$ (for $1 \leq k \leq {\cal N}_j$), and
use the $p$-triangle inequality, to obtain
\begin{equation}
\label{sjsumalp}
|s_j(k)|^p \leq |u^0_j(k)|^p +  
\sum_{\sigma_0 \leq \sigma \leq jc}\: \sum_{l(Q)=2^{-\sigma},3Q \subseteq W} |s_Q|^p |u^Q_j(k)|^p
+ \sum_{\sigma > jc}\: \sum_{l(Q)=2^{-\sigma},3Q \subseteq W} |s_Q|^p |u^Q_j(k)|^p.
\end{equation}
Let $A_j$, $B_{\sigma}$ be as above.
Integrating both sides of (\ref{sjsumalp}) over $Y_j$, using (\ref{uqjest1}), (\ref{uqjest2}), (\ref{uqjest3}), and using
Lemma \ref{fjan4} (taking ``$p$'' in that lemma to be $1$), we find
\[ A_j^p \leq C({\cal T}^p a^{-2jlp} + \sum_{\sigma_0 \leq \sigma \leq jc} a^{-2jlp} 2^{2\sigma lp}B_{\sigma}^p
+ \sum_{\sigma > jc} a^{j(2l+n)p-jn} 2^{-\sigma(2l+n)p+\sigma n}B_{\sigma}^p). \]
Let $A^{\alpha}_j$, $B^{\alpha}_{\sigma}$ be as above.  We find that

\begin{equation}
\label{ajbnalpest}
(A_j^{\alpha})^p \leq C({\cal T}^p a^{-j(2l-\alpha)p}+ \sum_{\sigma=\sigma_0}^{\infty} K(j,\sigma) (B^{\alpha}_{\sigma})^p),
\end{equation}

where now

\[ K(j,\sigma) = a^{-j(2l-\alpha)p}2^{\sigma(2l-\alpha)p} \mbox{    if } \sigma_0 \leq \sigma \leq jc, \]

\[ K(j,\sigma) = a^{j[(2l+n+\alpha)p -n]}2^{-\sigma[(2l+n+\alpha)p -n]} \mbox{    if } \sigma_0 \geq jc. \]

By (\ref{lgqgint}), $2l$ is more than $\max(\alpha,n/p-n-\alpha)$.  
Recall also that $c = \log_2 a$.  Thus, by Proposition \ref{lqconv},
$\|(A_j^{\alpha})^p\|_{q/p} \leq C({\cal T}^q+\|(B_{\sigma}^{\alpha})^p\|_{q/p})$.
Upon raising both sides to the $1/q$ power, one again obtains (\ref{bescmp1}), as desired.\\
\ \\
{\bf Remark} Presumably the assumptions in the last lemma can be weakened: assuming only
$2m > n(1/p-1)_+ - \alpha$, we conjecture that one should be able to replace the assumption (\ref{var2lph})
by the assumption that $\varphi^j_k = (a^{-2j}\Delta)^m\Phi^j_k$, where $\Phi^j_k$ satisfies
(\ref{xphip}) for all $X \in {\mathcal P}^{2(l+m)}$ (not ${\mathcal P}^{4l}$).

\begin{lemma}
\label{besov2}
Fix $b > 0$.  Also
fix an integer $l \geq 1$ with
\begin{equation}
\label{lgqgintac}
2l > n(1/p-1)_+ - \alpha.
\end{equation}
where here $x_+ = \max(x,0)$.  Fix $M$ with $(M-n)p > n+1$ if $0 < p < 1$, $M-n > n+1$
otherwise.

If $0 < p < 1$, we also fix a number $\rho > 0$.
Then there exists $C > 0$ as follows. 

Say $j \in \ZZ$.  Select sets $E^j_k$ and points $x^j_k$ as in Lemma \ref{fjan3}.
If $0 < p < 1$, we assume that, for all $j,k$,
\begin{equation}
\label{ejkro}
\mu(E^j_k) \geq \rho a^{-jn}
\end{equation}
Suppose that, for each $j \geq 0$, and each $k$, 
$\varphi^j_k = (a^{-2j}\Delta)^l\Phi^j_k$, where $\Phi^j_k \in C^{\infty}({\bf M})$ satisfies the following conditions:
\[ |X\Phi^j_k(y)| \leq a^{j (\deg X + n)}(1+a^{j}d(y,x^j_k))^{-M} \mbox{  whenever } X \in {\mathcal P}^{4l}. \]
Suppose that $\{s_{j,k}: j \geq 0, 1 \leq k \leq {\cal N}_j\}$ satisfies 
\[ (\sum_{j = 0}^{\infty} a^{j\alpha q} [\sum_k \mu(E^j_k)|s_{j,k}|^p]^{q/p})^{1/q} < \infty. \]
Then $\sum_{j=0}^{\infty} \sum_k \mu(E^j_k) s_{j,k} \varphi^j_k$ converges in $B_p^{\alpha q}({\bf M})$, and
\begin{equation}
\label{besov2way}
\|\sum_{j=0}^{\infty} \sum_k \mu(E^j_k) s_{j,k} \varphi^j_k\|_{B_p^{\alpha q}} \leq
C(\sum_{j = 0}^{\infty} a^{j\alpha q} [\sum_k \mu(E^j_k)|s_{j,k}|^p]^{q/p})^{1/q}.
\end{equation}
\end{lemma}
{\bf Proof} In \cite{SS}, Seeger and Sogge gave an equivalent characterization of $B_p^{\alpha q}$. 
(We change their notation a little; what we shall call $\beta_{k-1}(s^2)$, they called $\beta_k(s)$.)
Choose $\beta_0 \in C_c^{\infty}((1/4,16))$, with the property that
for any $s > 0$, $\sum_{\nu=-\infty}^{\infty} \beta^2_0(2^{-2\nu}s) = 1$.  For $\nu \geq 1$, define
$\beta_{\nu} \in C_c^{\infty}((2^{2\nu-2},2^{2\nu+4}))$, by $\beta_{\nu}(s) = \beta_0(2^{-2\nu}s)$.
Also, for $s > 0$, define 
the smooth function $\beta_{-1}(s)$ by
$\beta_{-1}(s) = \sum_{\nu=-\infty}^{-1} \beta(2^{-2\nu}s)$.  (Note that $\beta_{-1}(s) = 0$
for $s > 4$.)
Then (\cite{SS}), for $F \in C^{\infty}({\bf M})$, $\|F\|_{B_p^{\alpha q}}$ is equivalent to 
the $l^q$ norm of the sequence $\{ 2^{\nu\alpha}\|\beta_{\nu}(\Delta)F\|_p: -1 \leq \nu \leq \infty\}$.  

For $\nu \geq -1$, let $J_{\nu}$ be the kernel of $\beta_{\nu}(\Delta)$.  Using the eigenfunction 
expansion of $\beta_{\nu}(\Delta)$ (see (\ref{kerexp}) of \cite{gmcw}), one
sees at once that $J_{-1}(x,y)$ is smooth in $(x,y)$.  Moreover, for $\nu \geq 0$, $\beta_{\nu} \in {\cal S}(\RR^+)$
and $\beta_{\nu}(0) = 0$; so $J_{\nu}$ is smooth as well.
For any integer $I \geq 0$, we set
$\beta_0^I(s) = \beta_0(s)/s^I$, and define $\beta_{\nu}^I$ by 
$\beta_{\nu}^I(s) = \beta_0^I(2^{-2\nu}s)$.  Then $\beta_{\nu}^I(s) = 2^{2\nu I}\beta_{\nu}(s)/s^I$, so that
$\beta_{\nu}(\Delta) = 2^{-2\nu I} \beta_{\nu}^I(\Delta)\Delta^I $.  Thus, if $J_{\nu}^I$ is the kernel of $\beta_{\nu}^I(\Delta)$,
we have 
\begin{equation}
\label{jjdyi}
J_{\nu}(x,y) = 2^{-2{\nu}I}\Delta_y^I J_{\nu}^I(x,y),
\end{equation}
where $\Delta_y$ means $\Delta$ as applied in the $y$ variable.  Also, by Lemma \ref{manmol},
since $\beta_{\nu}^I(\Delta) = \beta_0^I(2^{-2\nu}\Delta)$, we know the following:
for every pair of
$C^{\infty}$ differential operators $X$ (in $x$) and $Y$ (in $y$) on ${\bf M}$,
and for every integer $N \geq 0$, and for any fixed $I$, there exists $C$ such that for all $\nu$,
\begin{equation}
\label{jliest}
|XY J_{\nu}^I(x,y)| \leq C 2^{\nu(n+\deg X + \deg Y)}(1+2^\nu d(x,y))^{-N}
\end{equation}
for all $x,y \in {\bf M}$.

In proving (\ref{besov2way}), we may assume that for all but finitely many $j$, all $s_{j,k} = 0$.  For if we can
prove the inequality (\ref{besov2way}) in that case, it will follow at once that the partial sums of
$\sum_{j=0}^{\infty} \sum_k \mu(E^j_k) s_{j,k} \varphi^j_k$ form a Cauchy sequence in the quasi-Banach space
$B_p^{\alpha q}({\bf M})$.  Thus the series will converge in that quasi-Banach space, and moreover the inequality
(\ref{besov2way}) will hold in full generality.  In fact (\ref{besov2way}) shows that the convergence is unconditional.

With this assumption, we may let $F = \sum_{j=0}^{\infty} \sum_k \mu(E^j_k) s_{j,k} \varphi^j_k$.

Let $c = \log_a 2$.  For $\nu \geq -1$, we have

\begin{equation}
\label{beftot}
\beta_{\nu}(\Delta)F  = 
\sum_{0 \leq j \leq \nu c}\: \sum_k s_{j,k} \mu(E^j_k)\beta_{\nu}(\Delta)\varphi^j_k 
+ \sum_{j > \nu c}\: \sum_k s_{j,k} \mu(E^j_k)\beta_{\nu}(\Delta)\varphi^j_k. 
\end{equation}

Of course, in each term, $[\beta_{\nu}(\Delta)\varphi^j_k](z) = \int J_{\nu}(x,y) \varphi^j_k(y) d\mu(y)$.

Suppose that $x \in {\bf M}$.

Now say that $j > \nu c$, i.e., that $a^j \geq 2^{\nu}$.  Then by Lemma \ref{fjan2} 
(replacing ``$\sigma$'' in that lemma by $j$, ``$x_0$'' by $x^j_k$, replacing
``$\nu$'' in that lemma by $\nu c$,
and with $\Phi^j_k = a^{jn} \Phi$, $\varphi^j_k = a^{jn-2jl}\varphi_1$, and
$J_{\nu}(x,y) = 2^{\nu (n+2l)}\varphi_2(y)$, then
\begin{equation}
\label{bdfest1}
|\mu(E^j_k)\beta_{\nu}(\Delta)\varphi^j_k(x)| 
\leq Ca^{-2jl} 2^{\nu(n+2l)}\mu(E^j_k)(1+2^{\nu}d(x,x^j_k))^{n-M}.
\end{equation}

Say instead $0 \leq j \leq \nu c$, so that $a^j \leq 2^{\nu}$.  Select $I$ with
\[ 2I > \max(\alpha,2l). \]
Then by Lemma \ref{fjan2}
(replacing ``$\sigma$'' in that lemma by 
$\nu c$, replacing ``$\nu$'' in that lemma by $j$, 
and with $J^I_{\nu}(x,y) = 2^{\nu n}\Phi(y)$, $J_{\nu}(x,y) = 2^{\nu n -2\nu I}\varphi_1(y)$, 
and $\varphi^j_k = a^{j(n+2l)}\varphi_2$)
we have
\begin{equation}
\label{bdfest2}
|\mu(E^j_k)\beta_{\nu}(\Delta)\varphi^j_k(x)| 
\leq Ca^{j(n+2I)} 2^{-2\nu I}\mu(E^j_k)(1+a^j d(x,x^j_k))^{n-M},
\end{equation}
since $a^{2jl} \leq a^{2jI}$.

Say now $p \geq 1$.  From (\ref{beftot}) we obtain
\[ \|\beta_{\nu}(\Delta)F\|_p \leq 
\sum_{0 \leq j \leq \nu c}\: \|\sum_k s_{j,k} \mu(E^j_k)\beta_{\nu}(\Delta)\varphi^j_k\|_p
+ \sum_{j > \nu c}\: \|\sum_k s_{j,k} \mu(E^j_k)\beta_{\nu}(\Delta)\varphi^j_k\|_p. \] 

Let 
\[ A_{\nu} = \|\beta_{\nu}(\Delta)F\|_p,\:\:
B_{j} = [\sum_k \mu(E^j_k)|s_{j,k}|^p]^{1/p}. \]

Then by (\ref{bdfest1}), (\ref{bdfest2}), and Lemma \ref{fjan3}, we see that

\begin{equation}
\label{ajbnuest3}
A_{\nu} \leq C(\sum_{0 \leq j \leq \nu c} a^{2jI} 2^{-2\nu I}B_{j}
+ \sum_{j > \nu c} a^{-2jl} 2^{2\nu l}B_{j}).
\end{equation}

Now also write

\[ A^{\alpha}_{\nu} = 2^{\nu \alpha} A_{\nu},\:\: B^{\alpha}_{j} = a^{j \alpha} B_{j}. \]

Then, by (\ref{ajbnuest3}),

\begin{equation}
\label{ajbnalest3}
A_{\nu}^{\alpha} \leq C(\sum_{j=0}^{\infty} K(\nu,j) B^{\alpha}_{j}),
\end{equation}

where 

\[ K(\nu,j) = 2^{-\nu(2I-\alpha)}a^{j(2I-\alpha)} \mbox{    if } 0 \leq j \leq \nu c, \]

\[ K(\nu,j) = 2^{\nu(2l+\alpha)}a^{-j(2l+\alpha)} \mbox{    if } j > \nu c. \]

By (\ref{lgqgintac}), $2l$ is more than $-\alpha$, and we have also taken $I$ to satisfy $2I > \alpha$. 
Recall also that $c = \log_a 2$.  Thus, by Proposition \ref{lqconv},
$\|(A_{\nu}^{\alpha})\|_q \leq C\|(B_{j}^{\alpha})\|_q$.  Consequently the $l^q$ norm of the
sequence $\{2^{\nu \alpha} \|\beta_{\nu}(\Delta)F\|_p\}$ is less than or equal to 
$C(\sum_{j = 0}^{\infty} a^{j\alpha q} [\sum_k \mu(E^j_k)|s_{j,k}|^p]^{q/p})^{1/q}$, which, because
of the result of Seeger-Sogge, gives the lemma, at least if $p \geq 1$.

If instead $0 < p < 1$, we evaluate each side of (\ref{beftot}) at $x$ (for each $x \in {\bf M}$), and
use the $p$-triangle inequality, to obtain
\[ |\beta_{\nu}(\Delta)F(x)|^p  \leq
\sum_{0 \leq j \leq \nu c}\: \sum_k |s_{j,k}|^p \mu(E^j_k)^p|\beta_{\nu}(\Delta)\varphi^j_k(x)|^p 
+ \sum_{j > \nu c}\: \sum_k |s_{j,k}|^p \mu(E^j_k)^p|\beta_{\nu}(\Delta)\varphi^j_k(x)|^p \]
so, by the assumption (\ref{ejkro}),
\begin{align}
&|\beta_{\nu}(\Delta)F(x)|^p \nonumber \\
&\leq  C (\sum_{0 \leq j \leq \nu c}\: \sum_k |s_{j,k}|^p a^{jn(1-p)}\mu(E^j_k)|\beta_{\nu}(\Delta)\varphi^j_k(x)|^p 
+ \sum_{j > \nu c}\: \sum_k |s_{j,k}|^p a^{jn(1-p)}\mu(E^j_k)|\beta_{\nu}(\Delta)\varphi^j_k(x)|^p) \label{beftotp}
\end{align}

Let $A_{\nu}$, $B_{j}$ be as above.
Integrating both sides of (\ref{beftotp}) over ${\bf M}$, using (\ref{bdfest1}), (\ref{bdfest2}), and Lemma \ref{fjan3}
(taking ``$p$'' in that lemma to be $1$), we find
\begin{equation}
\label{anupbjp}
A_{\nu}^p \leq C(\sum_{0 \leq j \leq \nu c} a^{2jIp} 2^{-2\nu Ip}B_{j}^p
+ \sum_{j > \nu c} a^{-j(2l+n)p+jn} 2^{\nu(2l+n)p-\nu n}B_{j}^p).
\end{equation}
Let $A^{\alpha}_{\nu}$, $B^{\alpha}_{j}$ be as above.  We find that

\begin{equation}
\label{ajbnalpest2}
(A_{\nu}^{\alpha})^p \leq C(\sum_{j=0}^{\infty} K(\nu,j) (B^{\alpha}_{j})^p),
\end{equation}

where now

\[ K(\nu,j) = 2^{-\nu(2I-\alpha)p}a^{j(2I-\alpha)p} \mbox{    if } 0 \leq \sigma \leq \nu c, \]

\[ K(\nu,j) = 2^{\nu[(2l+n+\alpha)p -n]}a^{-j[(2l+n+\alpha)p -n]} \mbox{    if } j \geq \nu c. \]

By (\ref{lgqgintac}), $2l$ is more than $n/p-n-\alpha$; also $2I > \alpha$.
Recall also that $c = \log_2 a$.  Thus, by Proposition \ref{lqconv},
$\|(A_{\nu}^{\alpha})^p\|_{q/p} \leq C(\|(B_{j}^{\alpha})^p\|_{q/p})$.
Upon raising both sides to the $1/q$ power, one again obtains (\ref{besov2way}), as desired.\\

For any $x \in {\bf M}$, and any integer $I, J \geq 1$, we let
\begin{equation}
\label{mxdefL}
{\mathcal M}_{x,t}^{IJ} = \{\varphi \in C_0^{\infty}({\bf M}):  
t^{n+\deg Y} |(\frac{d(x,y)}{t})^N Y\varphi(y)|\leq 1 
\mbox{ whenever } y \in {\bf M}, \: 0 \leq N \leq J \mbox{ and } 
Y \in {\mathcal P}^I\}.
\end{equation}
(This space is a variant of 
a space of molecules, as defined earlier in \cite{gil} and \cite{han2}.
In the notation of our earlier article \cite{gm2}, ${\mathcal M}_{x,t}^{n+2,1} = {\mathcal M}_{x,t}$.)
Note that, if for each $x \in {\bf M}$, we 
define the functions $\varphi^t_x, \psi^t_x$ on ${\bf M}$ by 
$\varphi^t_x(y) = K_t(x,y)$ and $\psi^t_x(y) = K_t(y,x)$ (notation as in 
Lemma \ref{manmol}), then by Lemma \ref{manmol}, for each
$I, J \geq 1$, 
there exists $C_{IJ} > 0$ such that $\varphi^t_x$ and $\psi^t_x$ are in $C_{IJ}{\mathcal M}_{x,t}^{IJ}$
{\em for all} $x \in {\bf M}$ and all $t > 0$.\\

We recall Theorem \ref{sumopthmfr} of our earlier article \cite{gm2}:
\begin{theorem}
\label{sumopthm}
Fix $a>1$.  Then there exists $C_1, C_2 > 0$ as follows.

For each $j \in \ZZ$, write ${\bf M}$ as a finite disjoint union of 
measurable subsets $\{E^j_k: 1 \leq k \leq {\cal N}_j\}$, each of diameter less than $a^{-j}$.
For each $j,k$, select any $x^j_k \in E^j_k$, and select $\varphi^j_k$,
 $\psi^j_k$ with \\$\varphi^j_k , \psi^j_k\in \mathcal{M}^{n+2,1}_{x^j_k, a^{-j}}$.
For $F \in C^1({\bf M})$, we claim that we may define 
\begin{equation}
\label{sumopdf2}
SF = S_{\{\varphi^j_k\},\{\psi^j_k\}}F = 
\sum_{j}\sum_k \mu(E^j_k) \langle  F, \varphi^j_k  \rangle  \psi^j_k. 
\end{equation}
$($Here, and in similar equations below, the sum in $k$ runs from $k = 1$ to $k = {\cal N}_j$.$)$
Indeed: 
\begin{itemize}
\item[$(a)$] For any $F\in C^1({\bf M})$, the series defining $SF$ converges absolutely, uniformly on $\bf{M}$,
\item[$(b)$] $\parallel SF\parallel_2\leq C_2\parallel F\parallel_2$ for all
 $F\in C^1(\bf{M})$.\\
Consequently, $S$ extends to be a bounded operator on $L^2({\bf M})$, with norm less than
or equal to $C_2$.
\item[$(c)$] If $F \in L^2({\bf M})$, then
\begin{equation}
\label{uncndl2}
SF = \sum_{j}\sum_k \mu(E^j_k) \langle  F, \varphi^j_k  \rangle  \psi^j_k
\end{equation}
where the series converges unconditionally. 
\item[$(d)$] If $F,G \in L^2({\bf M})$, then 
\begin{equation}\label{eq:wk}
 \langle  SF,G  \rangle = 
\sum_{j}\sum_{k}\mu(E^j_k) \langle  F, \varphi^j_k \rangle   \langle \psi^j_k,G \rangle ,
\end{equation}
where the series converges absolutely.\\
\end{itemize}
\end{theorem}

This result, which was proved by use of the $T(1)$ theorem, explains the $L^2$ theory of the summation 
operator $S$.  On Besov spaces, we have the following result for summation operators, where now we consider the sum over
nonnegative $j$:

\begin{theorem}
\label{sumopbes}
Fix an integer $l \geq 1$ with $2l > \max(n(1/p-1)_+ - \alpha, \alpha)$.  Also fix $J$
with $(J-2l-n)p > n+1$ if $0 < p < 1$, $J-2l-n > n+1$ if $p \geq 1$.

Then there exists $C > 0$ as follows.
For each integer $j \geq 0$, write ${\bf M}$ as a finite disjoint union of 
measurable subsets $\{E^j_k: 1 \leq k \leq {\cal N}_j\}$, each of diameter less than $a^{-j}$,
select $x^j_k \in E^j_k$, select 
$\Phi^j_k, \Psi^j_k \in \mathcal{M}_{x^j_k, a^{-j}}^{4l,J}$ and set
$\varphi^j_k = (a^{-2j}\Delta)^l\Phi^j_k$, $\psi^j_k = (a^{-2j}\Delta)^l\Psi^j_k$. 

By Theorem \ref{sumopthm}, we may then define 
\begin{equation}
\label{sumopdf4}
S^{\prime}F = S^{\prime}_{\{\varphi^j_k\},\{\psi^j_k\}}F = 
\sum_{j=0}^{\infty}\sum_k \mu(E^j_k) \langle  F, \varphi^j_k  \rangle  \psi^j_k,
\end{equation}
at first for $F \in C^1({\bf M})$, and Theorem \ref{sumopthm} applies.\\

Then, if $F \in B_p^{\alpha q}$, the series in (\ref{sumopdf4}) converges in $B_p^{\alpha q}$, to 
a distribution $S'F \in B_p^{\alpha q}$, such that $(S'F)(1) = 0$.  Moreover, $\|S'F\|_{B_p^{\alpha q}}
\leq C\|F\|_{B_p^{\alpha q}}$.
\end{theorem}
{\bf Proof}  Setting $s_{j,k} = \langle F, \varphi^j_k \rangle$, we see by Lemmas \ref{besov1} and 
\ref{besov2} that the series in (\ref{sumopdf4}) converges to a distribution $S'F$ in $B_p^{\alpha q}$,
and that 
\begin{equation}
\label{sumbes1}
\|S'F\|_{B_p^{\alpha q}}  = \|\sum_{j=0}^{\infty}\sum_k \mu(E^j_k) \langle F, \varphi^j_k \rangle  \psi^j_k\|_{B_p^{\alpha q}} 
\leq C(\sum_{j = 0}^{\infty} a^{j\alpha q} 
[\sum_k \mu(E^j_k)|\langle F, \varphi^j_k \rangle|^p]^{q/p})^{1/q} \leq C\|F\|_{B_p^{\alpha q}},
\end{equation}
{\em provided} that $p \geq 1$, or alternatively that $0 < p < 1$ and (\ref{ejkro}) holds.  If $0 < p < 1$ and
(\ref{ejkro}) does not hold, we at least know that the second inequality in (\ref{sumbes1}) holds.  For the
first inequality, we need to regroup.

Say then that $0 < p < 1$.  
By the discussion before Theorem \ref{framainfr} of \cite{gm2}, for each $j \geq 0$,
there exists a finite covering of ${\bf M}$ by 
disjoint measurable sets ${\cal F}^j_1,\ldots,{\cal F}^j_{L(j)}$, such that whenever $1 \leq i \leq L(j)$, there
is a $y^j_i \in {\bf M}$ with $B(y^j_i,2a^{-j}) \subseteq {\cal F}^j_i \subseteq B(y^j_i,4a^{-j})$.  
Fix $j$ for now.  For each $k$ with $1 \leq k \leq {\cal N}_j$, select a number $i_{j,k}$ with $1 \leq i_{j,k} \leq L(j)$,
such that ${\cal F}^j_{i_{j,k}} \cap E^j_k \neq \oslash$.  For $1 \leq i \leq L(j)$, let $S_{j,i} = \{k: i_{j,k} = i\}$,
and then let
\[ {\cal E}^j_i = \cup_{k \in S_{j,i}} E^j_k.  \]
Then we have that 
\[ B(y^j_i,a^{-j}) \subseteq {\cal E}^j_i \subseteq B(y^j_i,5a^{-j}).  \]
The second inclusion here is evident from the facts that, firstly, ${\cal F}^j_i \subseteq B(y^j_i,4a^{-j})$, 
secondly, that ${\cal F}^j_{i} \cap E^j_k \neq \oslash$ whenever $k \in S_{j,i}$, and thirdly, that each
$E^j_k$ has diameter less than $a^{-j}$.  For the first inclusion, say $x \in B(y^j_i,a^{-j})$; we need
to show that $x \in {\cal E}^j_i$.  Choose $k$ with $x \in E^j_k$.  Since $E^j_k$ has diameter less than $a^{-j}$,
surely $E^j_k \subseteq B(y^j_i,2a^{-j}) \subseteq {\cal F}^j_i$.  Thus ${\cal F}^j_i$ is the only
one of the sets ${\cal F}^j_1,\ldots,{\cal F}^j_{L(j)}$ which $E^j_k$ intersects, and so $k$ must be in 
$S_{j,i}$.  Accordingly $x \in E^j_k \subseteq {\cal E}^j_i$, as claimed.

For each $j, i$ let 
\[ r^{j,i} = \max_{k \in S_{j,i}} |\langle F, \varphi^j_k \rangle|. \]

Then, for some $k \in S_{j,i}$, $r^{j,i} = |\langle F, \varphi^j_k \rangle|$.
Now $d(x^j_k, y^j_i) \leq 5a^{-j}$, so it is evident from (\ref{mxdefL}) that, for some
absolute constant $C_0$, $\mathcal{M}_{x^j_k, a^{-j}}^{4l,J} \subseteq C_0\mathcal{M}_{y^j_i, a^{-j}}^{4l,J}$;
in particular, $\Phi^j_k \in C_0\mathcal{M}_{y^j_i, a^{-j}}^{4l,J}$.
Also, $\mbox{diam}{\cal E}^j_i \leq 10a^{-j}$. Thus, by Lemma \ref{besov1}, 
\begin{equation}
\label{rjifbpq}
(\sum_{j = 0}^{\infty} a^{j\alpha q} [\sum_i \mu({\cal E}^j_i)|r^{j,i}|^p]^{q/p})^{1/q} \leq C\|F\|_{B_p^{\alpha q}}.
\end{equation}

Let ${\bf F}$ be any finite subset of $\{(j,k): j \geq 0, 1 \leq k \leq {\cal N}_j\}$.  Define 
$s_{j,k} = \langle F, \varphi^j_k \rangle$ if $(j,k) \in {\bf F}$; otherwise, let $s_{j,k} = 0$.
Also, for each $j,i$ with $1 \leq i \leq L(j)$, let
\[ s^{j,i} = \max_{k \in S_{j,i}} |s_{j,k}|. \]

Then only finitely many of the $s^{j,i}$ are nonzero, and we always have $0 \leq s^{j,i} \leq r^{j,i}$.
Therefore by (\ref{rjifbpq}),
in order to show the convergence of the series for $S'F$ in  $B_p^{\alpha q}$, and that
$\|S'F\|_{B_p^{\alpha q}} \leq C\|F\|_{B_p^{\alpha q}}$, it is sufficient to show that
\begin{equation}
\label{sjksjiok}
\|\sum_{j=0}^{\infty}\sum_k s_{j,k} \mu(E^j_k) \psi^j_k\|_{B_p^{\alpha q}} \leq
C(\sum_{j = 0}^{\infty} a^{j\alpha q} [\sum_i \mu({\cal E}^j_i)|s^{j,i}|^p]^{q/p})^{1/q},
\end{equation}
where here $C$ is independent of ${\bf F}$ (and our choices of $E^j_k, x^j_k, \Phi^j_k,
\Psi^j_k, {\cal E}^j_i$).

Let $G = \sum_{j=0}^{\infty}\sum_k  s_{j,k} \mu(E^j_k) \psi^j_k$.
Let $\beta_{\nu}$ and $c$ be as in the proof of Lemma \ref{besov2}.
We have from (\ref{beftot}) that, for each $x \in {\bf M}$, 
\begin{equation}
\label{beftotabs}
|\beta_{\nu}(\Delta)G(x)| \leq 
\sum_{0 \leq j \leq \nu c}\: \sum_{i=1}^{L(j)}\sum_{k \in S_{j,i}} \mu(E^j_k)|s_{j,k}| |\beta_{\nu}(\Delta)\psi^j_k(x)|
+ \sum_{j > \nu c}\: \sum_{i=1}^{L(j)}\sum_{k \in S_{j,i}}\mu(E^j_k)|s_{j,k}| \beta_{\nu}(\Delta)\psi^j_k(x)|.
\end{equation}
For each $j, i$, and each $x \in {\bf M}$, let
\[ H^{j,i}(x) = \max_{k \in S_{j,i}}|\beta_{\nu}(\Delta)\psi^j_k(x)|. \]

From (\ref{beftotabs}), we see that 
\begin{equation}
\label{beftotamlg}
|\beta_{\nu}(\Delta)G(x)| \leq 
\sum_{0 \leq j \leq \nu c}\: \sum_{i=1}^{L(j)} s^{j,i} \mu({\cal E}^j_i)H^{j,i}(x)
+ \sum_{j > \nu c}\: \sum_{i=1}^{L(j)}s^{j,i} \mu({\cal E}^j_i)H^{j,i}(x).
\end{equation}

Note once again that if $k \in S_{j,i}$, then $d(x^j_k, y^j_i) \leq 5a^{-j}$; and therefore,
$\Psi^j_k \in C_0\mathcal{M}_{y^j_i, a^{-j}}^{4l,J}$ for some absolute constant $C_0$.  Note also that
$H^{j,i}(x) = |\beta_{\nu}(\Delta)\psi^j_k(x)|$ for some $k \in S_{j,i}$.  Thus, the reasoning leading to
(\ref{bdfest1}) and (\ref{bdfest2}) shows that if $j > \nu c$, then 
\begin{equation}
\label{bdfest3}
\mu({\cal E}^j_i) H^{j,i}(x)
\leq Ca^{-2jl} 2^{\nu(n+2l)}\mu({\cal E}^j_i)(1+2^{\nu}d(x,y^j_i))^{n-J},
\end{equation}
while if we select $I$ with $2I > \max(\alpha, 2l)$, and if $0 \leq j \leq \nu c$, then
\begin{equation}
\label{bdfest4}
\mu({\cal E}^j_i) H^{j,i}(x)
\leq Ca^{j(n+2I)} 2^{-\nu(2I)}\mu({\cal E}^j_i)(1+a^j d(x,y^j_i))^{n-J}.
\end{equation}

Since $B(y^j_i,a^{-j}) \subseteq {\cal E}^j_i$ for all $j, i$, there exists $\rho > 0$ such that
$\mu({\cal E}^j_i) \geq \rho a^{-jn}$ for all $j \geq 0$.  Consequently, the reasoning leading to (\ref{beftotp})
shows that 

\[ |\beta_{\nu}(\Delta)G(x)|^p 
\leq  C (\sum_{0 \leq j \leq \nu c}\: \sum_i |s^{j,i}|^p a^{jn(1-p)}\mu({\cal E}^j_i) |H^{j,i}(x)|^p
+ \sum_{j > \nu c}\: \sum_k |s^{j,i}|^p a^{jn(1-p)}\mu({\cal E}^j_i) |H^{j,i}(x)|^p). \]

Now set
\[ A_{\nu} = \|\beta_{\nu}(\Delta)G\|_p,\:\:
B_{j} = [\sum_i \mu({\cal E}^j_i)|s^{j,i}|^p]^{1/p},\:\:
A^{\alpha}_{\nu} = 2^{\nu \alpha} A_{\nu},\:\: B^{\alpha}_{j} = a^{j \alpha} B_{j}. \]

Just as in the proof of Lemma \ref{besov2}, 
(\ref{bdfest3}), (\ref{bdfest4}) and Lemma \ref{fjan3} (taking ''$p$'' in that lemma to be $1$) show 
that (\ref{anupbjp}) and (\ref{ajbnalpest2}) hold, with $K(\nu,j)$ as just after (\ref{ajbnalpest2}).
This again gives $\|(A_{\nu}^{\alpha})^p\|_{q/p} \leq C(\|(B_{j}^{\alpha})^p\|_{q/p})$.
Upon raising both sides to the $1/q$ power, one obtains (\ref{sjksjiok}), as claimed.

Finally, the fact that $(S'F)(1) = 0$ follows from the fact that each term of the summation in
(\ref{sumopdf4}) vanishes when applied to $1$ (since the assumption that $l \geq 1$ implies that
each $\psi^j_k$ has integral zero),
and the fact that the series in (\ref{sumopdf4}) converges in $B_p^{\alpha q}({\bf M})$, and hence in 
${\cal E}'({\bf M})$.  This completes the proof.

\begin{definition}
\label{bes0}
We let $B_{p,0}^{\alpha q}({\bf M}) = \{F \in B_p^{\alpha q}({\bf M}): F1 = 0\}$.  (Here $F1$ is the
result of applying the distribution $F$ to the constant function $1$.)
\end{definition}

\begin{theorem}
\label{besmain1}
Say $c_0 , \delta_0 > 0$.  Fix an integer $l \geq 1$ with $2l > \max(n(1/p-1)_+ - \alpha, \alpha)$.
Say $f_0 \in {\cal S}(\RR^+)$, and let $f(s) = s^l f_0(s)$.  
Suppose also that the Daubechies condition (\ref{daub}) holds.
Then there exist constants $C > 0$  and $0 < b_0 < 1$ as follows:\\

Say $0 < b < b_0$.
Suppose that, for each $j$, we can write ${\bf M}$ as a finite disjoint union of measurable sets 
$\{E_{j,k}: 1 \leq k \leq N_j\}$, and that (\ref{diamleq}), (\ref{measgeq}) hold.

(In the notation of Theorem \ref{framainfr} of \cite{gm2}, this is surely possible if
$c_0 \leq c_0^{\prime}$ and $\delta_0 \leq 2\delta$.)
Select $x_{j,k} \in E_{j,k}$ for each $j,k$.   
For $t > 0$, let $K_t$ be the kernel of $f(t^2\Delta)$.  Set
\[ \varphi_{j,k}(y) =  \overline{K}_{a^j}(x_{j,k},y). \]

By Lemma (\ref{manmol}), there is a constant $C_0$ (independent of the
choice of $b$ or the $E_{j,k}$), such that 
$\varphi_{j,k} \in C{\mathcal M}_{x_{j,k},a^j}$
for all $j,k$.  Thus, we may form the summation operator $S$ with

\begin{equation}
\label{sumjopdfbes}
SF = S_{\{\varphi_{j,k}\},\{\varphi_{j,k}\}}F = 
\sum_{j}\sum_k \mu(E_{j,k}) \langle  F, \varphi_{j,k}  \rangle  \varphi_{j,k}, 
\end{equation}
at first for $F \in C^1({\bf M})$, and Theorem \ref{sumopthm} applies.

Then:\\
If $F \in B_{p,0}^{\alpha q}$, then the series in (\ref{sumjopdfbes}) converges in $B_{p,0}^{\alpha q}$,
and $S: B_{p,0}^{\alpha q} \rightarrow B_{p,0}^{\alpha q}$ is bounded and {\em invertible}.  We have
$\|SF\|_{B_p^{\alpha q}} \leq C\|F\|_{B_p^{\alpha q}}$ 
and $\|S^{-1}F\|_{B_p^{\alpha q}} \leq C\|F\|_{B_p^{\alpha q}}$ 
for all $F \in B_{p,0}^{\alpha q}$.
\end{theorem}
{\bf Proof} By the last sentence of Lemma \ref{manmol},
if we let $K^0_t(x,y)$ be the kernel of $f_0(t^2\Delta)$,
and if we let $\eta_{j,k}(y) = \overline{K}^0_{a^j}(x_{j,k},y)$ for $j \leq 0$, then for every $I,J$ there exists
$C_{IJ}$ with $\eta_{j,k} \in C_{IJ}{\mathcal M}^{IJ}_{x_{j,k},a^j}$.  Also (for instance, by 
looking at eigenfunction expansions, as in (\ref{kerexp}) of \cite{gmcw}), one has that
$\varphi_{j,k} = (a^{2j}\Delta)^l \eta_{j,k}$.  Thus, by Theorem \ref{sumopbes}, if $F \in B_p^{\alpha q}$, the series

\begin{equation}
\label{sumjopdfbesneg} 
\sum_{j \leq 0}\sum_k \mu(E_{j,k}) \langle  F, \varphi_{j,k}  \rangle  \varphi_{j,k} := S'F
\end{equation}
converges in $B_{p,0}^{\alpha q}$, and $\|S'F\|_{B_{p,0}^{\alpha q}} \leq C\|F\|_{B_p^{\alpha q}}$.  

We shall next show that, if $F \in B_p^{\alpha q}$, then the series

\begin{equation}
\label{sprprf} 
\sum_{j > 0}\sum_k \mu(E_{j,k}) \langle  F, \varphi_{j,k}  \rangle  \varphi_{j,k}
\end{equation}
also converges in $B_{p,0}^{\alpha q}$, and if we call the sum of this series $S''F$, then
$\|S''F\|_{B_{p,0}^{\alpha q}} \leq C\|F\|_{B_p^{\alpha q}}$.  

Since $S = S' + S''$,
this will tell us that $\|SF\|_{B_{p,0}^{\alpha q}} \leq C\|F\|_{B_p^{\alpha q}}$.

More generally, let us show the following:\\
\ \\
(*) There exist $N, C_0$ (independent of our choices of $b$,
$E_{j,k}$, $x_{j,k}$) as follows.  Suppose $\psi_{j,k}, \Psi_{j,k}$ are smooth functions
on ${\bf M}$ (for all $j > 0$ and $1 \leq k \leq N(j)$), which satisfy
\begin{equation}
\label{jposway}
\|\psi_{j,k}\|_{C^N} \leq C_1 a^{-2j},\:\: \|\Psi_{j,k}\|_{C^N} \leq C_1 a^{-2j},
\end{equation}
for all $j,k$.  Then, if $F \in B_p^{\alpha q}$,  the series
\begin{equation}
\label{jposway1}
\sum_{j > 0}\sum_k \mu(E_{j,k}) \langle  F, \psi_{j,k}  \rangle  \Psi_{j,k}
\end{equation}
converges in $B_{p,0}^{\alpha q}$, and if we call the sum of this series $S''F = 
S^{\prime \prime}_{\{\psi_{j,k}\},\{\Psi_{j,k}\}}F$, then 
\begin{equation}
\label{jposway2}
\|S''F\|_{B_{p,0}^{\alpha q} }\leq C_0 C_1\|F\|_{B_p^{\alpha q}}.
\end{equation}  

Note that our $\varphi_{j,k}$ do satisfy (\ref{jposway}), since, as we noted in section 4 of \cite{gmcw}, 
$\lim_{t \rightarrow \infty} t K_{\sqrt t}(x,y) = 0$ in $C^{\infty}({\bf M} \times {\bf M})$.\\

To prove (*), we need only note that, by (\ref{ccNbes}) and a partition of unity argument,
one has (for some $N$) a continuous inclusion $C^N({\bf M}) \subseteq B_p^{\alpha q}({\bf M})$.
Also, since the inclusion $B_p^{\alpha q} \subseteq {\cal S}'(\RR^n)$ is continuous, we have a
continuous inclusion $B_p^{\alpha q}({\bf M}) \subseteq {\cal E}'(\RR^n)$.  In particular, for some $N$, 
\begin{equation}
\label{bescn}
\|G\|_{B_p^{\alpha q}} \leq C\|G\|_{C^N}
\end{equation}
for all $G \in C^N$,
while 
\begin{equation}
\label{besepcn}
|\langle F, G \rangle| \leq C\|F\|_{B_p^{\alpha q}}\|G\|_{C^N}
\end{equation}
for all $F \in B_p^{\alpha q}$ and all
$G \in C^{\infty}$.  In particular, in (\ref{jposway1}), 
\[\sum_k \mu(E_{j,k}) |\langle  F, \psi_{j,k} \rangle|\|\Psi_{j,k}\|_{C^N}
\leq CC^2_1 a^{-4j}\|F\|_{B_p^{\alpha q}}\sum_k \mu(E_{j,k}) 
\leq CC^2_1 \mu({\bf M}) a^{-4j}\|F\|_{B_p^{\alpha q}}.\]
Therefore the series in (\ref{jposway1}) converges 
absolutely in $C^N$, hence in $B_p^{\alpha q}$, and we have the estimate 
(\ref{jposway2}) as well.\\
 
To complete the proof, we return to the notation we used in the proof of Theorem \ref{framainfr}
of \cite{gm2} (taking
``${\mathcal J}$'' in that proof to be $\ZZ$, and setting $Q = Q^{\ZZ}$).  We wish to show that
$Q$, when restricted to $C^{\infty}({\bf M})$, has a bounded extension to an operator
$Q: B_p^{\alpha q} \rightarrow B_{p,0}^{\alpha q}$, and that, as operators on $B_p^{\alpha q}$, 
\[ \|Q-S\| \leq Cb \]
(where as usual $C$ is independent of our choices of $b$, the $E_{j,k}$ and the $x_{j,k}$).
Now, $C^{\infty}$ is dense in $B_p^{\alpha q}$ (for instance, by Theorem 7.1 (a) of \cite{FJ1};
the constructions in that paper show that the building blocks can be taken to be smooth).  Thus
it is enough to show that 
$\|(Q-S)F\| \leq Cb\|F\|$ for all $F \in C^{\infty}$ (where, until further notice, $\|\:\|$ means
the $B_p^{\alpha q}$ norm).  

Say then that $F \in C^{\infty}$.  As in the proof of Theorem \ref{framainfr} of \cite{gm2}, we may assume that 
$\delta_0 \leq \delta$.   Again put $\Omega_b = \log_a(\delta_0/b)$.
In the proof of Theorem \ref{framainfr} of \cite{gm2}, we have obtained a formula
for $\langle (Q-S)F, F \rangle$.  This expression may be polarized, and a formula for $(Q-S)F$ may be 
obtained from it.  We see that $(Q-S)F = \sum_{i=1}^N I_i + II$, where in the notation of the proof of Theorem 
\ref{framainfr} of \cite{gm2}, 

\[ I_i = b \int_{\mathcal B}
\int_0^1 [ S_{\{\varphi_{j,k}^{w,s}\},\{\psi_{j,k}^{w,s}\}}F
+   S_{\{\psi_{j,k}^{w,s}\},\{\varphi_{j,k}^{w,s}\}}F ] dsdw \]

and where
\[ II = \sum_{j \geq \Omega_b} \int_{\bf M} \langle  F, \Phi_{x,a^j} \rangle \Phi_{\cdot,a^j} d\mu(x)
+ \sum_{j \geq \Omega_b} \sum_k \mu(E_{j,k}) \langle F, \varphi_{j,k} \rangle \varphi_{j,k}. \]

(Note that $I_i$ and $II$ are functions.  In the first term of $II$ we have represented the independent 
variable as $\cdot$.)

In $I_i$ we note, from the explicit expressions for the $\varphi_{j,k}^{w,s}$ and $\psi_{j,k}^{w,s}$, 
and by using $K^0_t$ as in the first paragraph of this proof,
that we can write $\varphi_{j,k}^{w,s} = (a^{-2j}\Delta)^l\eta_{j,k}^{w,s}$,
$\psi_{j,k}^{w,s} = (a^{-2j}\Delta)^l\Psi_{j,k}^{w,s}$, where for any $I,J$ there exists $C_{IJ}$ with
$\eta_{j,k}^{w,s}, \Psi_{j,k}^{w,s} \in C_{IJ}\mathcal{M}_{x_{j,k}, a^j}^{IJ}$.  By what we have already 
seen in this proof, this implies that $\|I_i\| \leq Cb\|F\|$.  Also, $I_i(1)=0$.

In $II$ we note that, for some $C$, $\|\Phi_{z,a^j}\|_{C^N} \leq Ca^{-2j}$ (for all $j \geq \Omega_b$ and 
all $z \in {\bf M}$), while $\|\varphi_{j,k}\|_{C^N} \leq Ca^{-2j}$ (for all $j \geq \Omega_b$ and all $k$).
By (\ref{bescn}) and (\ref{besepcn}), we see that 
\[ \|II\| \leq C \|F\| \sum_{j \geq \Omega_b}a^{-4j} \leq Cb \|F\|, \]
for $0 < b < 1$, since $\Omega_b = \log_a(\delta_0/b)$.  Moreover, $II(1) = 0$.

We have, then, that for $F \in C^{\infty}$, $\|QF\| \leq \|SF\| + Cb\|F\|$.  So 
$Q: B_p^{\alpha q} \rightarrow B_p^{\alpha q}$ is bounded.  
Also, we see that $\|Q-S\| \leq Cb$.   Further, since $Q=S+\sum_{i=1}^N I_i + II$, we know $Q(1)=0$.
To complete the proof, it suffices to show that 
$Q: B_{p,0}^{\alpha q} \rightarrow B_{p,0}^{\alpha q}$ is invertible.  Indeed, if $\|\:\|$ now
denotes the norm of an operator on $B_{p,0}^{\alpha q}$, we will then have that 
$\|I-Q^{-1}S\| \leq Cb\|Q^{-1}\|$, and the theorem will follow for $b$ sufficiently small.

For $\lambda \in \RR$, let $H(\lambda) = \sum_{j=-\infty}^{\infty} |f|^2(a^{2j} \lambda^2)
= \sum_{j=-\infty}^{\infty} (a^{j}\lambda)^{4l}g_0(a^j \lambda)$, if we write $g_0(\xi) = |f_0|^2(\xi^2)$.
Since $g_0$ and all its derivatives are bounded and decay rapidly
at $\infty$, $H$ is a smooth even function on $\RR^+ \setminus \{0\}$.  
By Daubechies' criterion, $1/H = G$, say, is a smooth even function on $\RR^+ \setminus \{0\}$.

Now, let $u$ equal either $G$ or $H$.  Note that $u(\lambda) = u(a^{-1}\lambda)$,
so that $u$ is actually a bounded smooth function on $\RR^+ \setminus \{0\}$.  Moreover, 
if $\lambda > 0$, we may choose an integer $m$ with $a^{m} \leq \lambda \leq a^{m+1}$. 
Since $u(\lambda) = u(a^{-m}\lambda)$, for any $k$ we have
\[  |u^{(k)}(\lambda)| = a^{-km} |u^{(k)}(a^{-m}\lambda)| \leq a^{k}\lambda^{-k}M, \]
where $M = \max_{1 \leq \lambda \leq a} |u^{(k)}(\lambda)|$.
This implies that $\|\lambda^k u^{(k)}(\lambda)\|_{\infty} < \infty$
for any $k$.  

Now choose an even function $v \in C_c^{\infty}(\RR)$ with $\supp v \subseteq (-\lambda_1,\lambda_1)$
and $v \equiv 1$ in a neighborhood of $0$.  (Here $\lambda_1$ is the smallest positive 
eigenvalue of $\Delta$).  Then the even function $u_1 := (1-v)u \in C^{\infty}(\RR)$ is in 
$S_1^0(\RR)$, so $u_1(\sqrt{\Delta}) \in OPS^0_{1,0}({\bf M})$.  Moreover, $u_1(\sqrt{\Delta}) =
u(\sqrt{\Delta})$ on $C_0^{\infty}$ (:= $(I-P)C^{\infty}({\bf M})$).  Recall that, here, $u= G$ or $H$;
set $G_1 = (1-v)G$, $H_1 = (1-v)H$.

Note further that $Q = H(\sqrt{\Delta}) = H_1(\sqrt{\Delta})$,
first when acting on finite linear combinations of non-constant
eigenfunctions of $\Delta$, then
on $C_0^{\infty}$, since such finite linear combinations are dense in $C_0^{\infty}$,
and the operators are bounded on $L^2$.  Moreover $G_1(\sqrt{\Delta})H_1(\sqrt{\Delta}) = 
H_1(\sqrt{\Delta})G_1(\sqrt{\Delta}) = (H_1G_1)(\sqrt{\Delta}) = I$.
first when acting on finite linear combinations of non-constant
eigenfunctions of $\Delta$, then on $C_0^{\infty}$, and finally on $B_{p,0}^{\alpha q}$, since by
\cite{P}, operators in $OPS^0_{1,0}({\bf M})$ are bounded on $B_p^{\alpha q}$.  This shows that
$Q$ is indeed invertible on $B_{p,0}^{\alpha q}$, and establishes the theorem.\\

In Theorem \ref{besmain1}, we have again required the condition (\ref{measgeq}), that
$\mu(E_{j,k}) \geq c_0(ba^j)^n$, whenever $ba^j < \delta_0$.
In the next theorem, if $0 < p < 1$, we will require a mild condition on $\mu(E_{j,k})$ if $ba^j \geq \delta_0$,
namely that
\begin{equation}
\label{measgeq2}
\mu(E_{j,k}) \geq {\cal C} \mbox{ whenever } ba^j \geq \delta_0,
\end{equation}
(for some ${\cal C} > 0$).  Sets $E_{j,k}$ which satisfy (\ref{diamleq}), (\ref{measgeq}) and
(\ref{measgeq2}) are easily constructed.   Indeed, set $t = ba^j/2 \geq \delta_0/2$, and,
as in the second bullet point prior to Theorem \ref{framainfr} of \cite{gm2}, select a 
finite covering of ${\bf M}$ by 
disjoint measurable sets $E_1,\ldots,E_N$, such that whenever $1 \leq k \leq N$, there
is a $y_k \in {\bf M}$ with $B(y_k,t) \subseteq E_k \subseteq B(y_k,2t)$.  Then, by (\ref{ballsn})
and (\ref{ballsn1}), there is a constant $c_0^{\prime}$, depending only on ${\bf M}$, such that 
$\mu(E_k) \geq c_0^{\prime}\min(\delta^n,(\delta_0/2)^n)$, as desired.


\begin{theorem}
\label{besmain2}
Say $c_0 , \delta_0, M_0, {\cal C} > 0$.  Fix an integer $l \geq 1$ with $2l > \max(n(1/p-1)_+ - \alpha, \alpha)$.
Say $f_0 \in {\cal S}(\RR^+)$, and let $f(s) = s^l f_0(s)$.  
Suppose also that the Daubechies condition (\ref{daub}) holds.
Then there exist constants $C_1 > 0$  and $0 < b_0 < 1$ as follows:\\

Say $0 < b < b_0$.  Then there exists a constant $C_2 > 0$ as follows:\\

Suppose that, for each $j$, we can write ${\bf M}$ as a finite disjoint union of measurable sets 
$\{E_{j,k}: 1 \leq k \leq N_j\}$, and that (\ref{diamleq}), (\ref{measgeq}) hold.  If $0 < p < 1$,
we suppose that (\ref{measgeq2}) holds as well.  Select $x_{j,k} \in E_{j,k}$ for each $j,k$. 
  
For $t > 0$, let $K_t$ be the kernel of $f(t^2\Delta)$.  Set
\[ \varphi_{j,k}(y) =  \overline{K}_{a^j}(x_{j,k},y). \]

Suppose $F$ is a distribution on ${\bf M}$ of order at most $M_0$, and that $F1=0$.  
Then the following are equivalent:\\
(i) $F \in B_{p,0}^{\alpha q}$;\\
(ii) $(\sum_{j = -\infty}^{\infty} a^{-j\alpha q} [\sum_k \mu(E_{j,k})
|\langle  F, \varphi_{j,k}  \rangle  |^p]^{q/p})^{1/q} < \infty$.  Further\\
\begin{equation}
\label{nrmeqiv1}
\|F\|_{B_p^{\alpha q}}/C_2 \leq (\sum_{j = -\infty}^{\infty} a^{-j\alpha q} [\sum_k \mu(E_{j,k})
|\langle  F, \varphi_{j,k}  \rangle  |^p]^{q/p})^{1/q} \leq C_1\|F\|_{B_p^{\alpha q}}.
\end{equation}
Moreover, if $p \geq 1$, then $C_2$ may be chosen to be independent of the choice of $b$
with $0 < b < b_0$.
\end{theorem}
{\bf Proof} Note first that, if $\eta_{j,k}$ is as in the first paragraph of the proof of Theorem 
\ref{besmain1}, then as we noted there for every $I,J$ there exists
$C_{IJ}$ with $\eta_{j,k} \in C_{IJ}{\mathcal M}^{IJ}_{x_{j,k},a^j}$, and 
$\varphi_{j,k} = (a^{2j}\Delta)^l \eta_{j,k}$. 

Say now that $F \in B_{p,0}^{\alpha q}$.  As we noted in the first paragraph of the proof of 
Theorem \ref{sumopbes}, it follows from Lemma \ref{besov1} that
\[ (\sum_{j = -\infty}^{0} a^{-j\alpha q} [\sum_k \mu(E_{j,k})
|\langle  F, \varphi_{j,k}  \rangle  |^p]^{q/p})^{1/q} \leq C\|F\|_{B_p^{\alpha q}}. \]

As for the terms in (ii) with $j > 0$, we note that, as in section 4 of \cite{gmcw}, 
$\lim_{t \rightarrow \infty} t^M K_{\sqrt t}(x,y) = 0$ in $C^{\infty}({\bf M} \times {\bf M})$,
for any $M$.  Consequently, for any $M,N$, 
\begin{equation}
\label{vphcnm}
\|\varphi_{j,k}\|_{C^N} \leq C a^{-Mj}.
\end{equation}

We choose any $M > -\alpha$. Then,
by (\ref{bescn}) and (\ref{besepcn}), we see that 

\[ (\sum_{j = 1}^{\infty} a^{-j\alpha q} [\sum_k \mu(E_{j,k})
|\langle  F, \varphi_{j,k}  \rangle  |^p]^{q/p})^{1/q} \leq
C\|F\|_{B_p^{\alpha q}} \mu(M)^{1/p} (\sum_{j = 1}^{\infty} a^{-j\alpha q - jMq})^{1/q} \leq
C\|F\|_{B_p^{\alpha q}}. \]

This proves that $(i) \Rightarrow (ii)$, and also establishes the rightmost inequality in 
(\ref{nrmeqiv1}).

Say, conversely, that $(ii)$ holds.  Our first step will be to show that the series for $SF$
(as in (\ref{sumjopdfbes}))
converges in $B_{p,0}^{\alpha q}$.  For this, 
it will be enough to show that the series for $S'F$ and $S''F$ (as in (\ref{sumjopdfbesneg})
and (\ref{sprprf}))each converge in $B_{p,0}^{\alpha q}$.  

Lemma \ref{besov2}, with $s_{j,k}$ in that Lemma being our $\langle F, \varphi_{-j,k} \rangle$,
implies that the series for $S'F$ does converge in $B_{p,0}^{\alpha q}$, and,
moreover, that
\begin{equation}
\label{sprfiii}
\|S'F\|_{B_p^{\alpha q}} \leq C(\sum_{j = -\infty}^{0} a^{-j\alpha q} [\sum_k \mu(E_{j,k})
|\langle  F, \varphi_{j,k}  \rangle  |^p]^{q/p})^{1/q}. 
\end{equation}

(If $0 < p < 1$, we need to note that, since we are assuming (\ref{measgeq}) and (\ref{measgeq2}), 
we have that (\ref{ejkro})
holds for some $\rho > 0$.  Since $\rho$ depends on $b$, the constant $C$ in (\ref{sprfiii}) depends on
$b$ as well, if $0 < p < 1$.)  

As for $S''F$, by (\ref{bescn}), it is enough to show that the series for it converges absolutely in $C^N$ (for any
fixed $N$).  By (\ref{vphcnm}), we need only show:\\
\ \\ 
(*)  Suppose that $\{r_{j,k}: 1 \leq j < \infty, 1 \leq k \leq N(j)\}$ are constants.  
If $M$ is sufficiently large, then 
\begin{equation}
\label{cnmjpq}
\sum_{j=1}^{\infty}a^{-Mj} \sum_k \mu(E_{j,k})|r_{j,k}|  \leq 
C(\sum_{j = 1}^{\infty} a^{-j\alpha q} [\sum_k \mu(E_{j,k}) |r_{j,k}|^p]^{q/p})^{1/q}.
\end{equation}
Here $C$ depends only on $c_0, \delta_0$ if $p \geq 1$ and only on $c_0, \delta_0, b$ if
$0 < p < 1$.\\
\ \\
To see this, let $a_j = \sum_k \mu(E_{j,k})|r_{j,k}|$, 
$d_j = [\sum_k \mu(E_{j,k}) |r_{j,k}|^p]^{1/p}$; we begin by showing that $a_j \leq C d_j$.
If $p \geq 1$, we let $Y_j$ be the finite measure space $\{1,\ldots,N_j\}$ with measure
$\lambda$, where $\lambda(\{k\}) = \mu(E_{j,k})$.  Applying H\"older's inequality on $Y_j$,
we find that $a_j \leq \mu({\bf M})^{1/p'} d_j$, as claimed.  If, instead, $0 < p < 1$, the 
$p$-triangle inequality implies that 
\[ a_j^p \leq \sum_k \mu(E_{j,k})^p|r_{j,k}|^p \leq C\sum_k \mu(E_{j,k}) |r_{j,k}|^p = Cd_j^p \]
where now $C$ depends on $b$.  (We have noted that, if $ba^j < \delta$ then 
$1 \leq \mu(E_{j,k})^{1-p}/[c_0 (ba^j)^n]^{1-p} \leq 
\mu(E_{j,k})^{1-p}/[c_0 b^n]^{1-p}$
by (\ref{measgeq}), while if
$ba^j \geq \delta$, then $1 \leq \mu(E_{j,k})^{1-p}/[{\cal C}]^{1-p}$ by (\ref{measgeq2}).)

To prove (\ref{cnmjpq}), it is enough, then, to show that, for $M$ sufficiently large,
\begin{equation}
\label{cnmjpq1}
\sum_{j=1}^{\infty}a^{-Mj} a_j  \leq C(\sum_{j = 1}^{\infty} a^{-j\alpha q} a_j^q)^{1/q},
\end{equation}
for the right side of (\ref{cnmjpq1}) is less than or equal to 
$C(\sum_{j = 1}^{\infty} a^{-j\alpha q} d_j^q)^{1/q}$, as we have seen, which gives (\ref{cnmjpq})
at once.  But (\ref{cnmjpq1}) is true for {\em any} nonnegative constants $a_j$ (for $M$ sufficiently
large), for the following reason.  If $0 < q \leq 1$, the $q$-triangle inequality tells us
that $(\sum_{j=1}^{\infty}a^{-j} a_j)^q \leq \sum_{j = 1}^{\infty} a^{-j\alpha q} a_j^q$, as claimed.
On the other hand, if $q > 1$, and $(-M + \alpha)q' < -1$, then (\ref{cnmjpq1}) follows by writing
$a^{-Mj} a_j = a^{(-M+\alpha)j}(a^{-\alpha j}a_j)$, and using H\"older's inequality.

In all, then, the series for $SF$
converges in $B_{p,0}^{\alpha q}$.  Moreover, if we call the sum of this series $S_0F$, we have from
(\ref{sprfiii}) and (\ref{cnmjpq}) that
\begin{equation}
\label{sofbpq}
\|S_0F\|_{B_p^{\alpha q}} \leq C(\sum_{j = -\infty}^{\infty} a^{-j\alpha q} [\sum_k \mu(E_{j,k})
|\langle  F, \varphi_{j,k}  \rangle  |^p]^{q/p})^{1/q},
\end{equation}
with $C$ independent of $b$ if $p \geq 1$.

To complete the proof we need only prove that if (ii) holds, and $b$ is sufficiently small,
then $F \in B_{p,0}^{\alpha q}$.  For then we will surely have that $S_0 F = SF$.  By Theorem \ref{besmain1},
$S: B_{p,0}^{\alpha q} \rightarrow B_{p,0}^{\alpha q}$ is invertible, with $\|S^{-1}\|$
independent of $b$ (for $b$ sufficiently small.)  Thus the
leftmost inequality in (\ref{nrmeqiv1}) will follow from (\ref{sofbpq}), with $C_2$ independent of $b$ if 
$p \geq 1$ (and $b$ is sufficiently small).

We have assumed that $F$ is a distribution of order at most $M_0$, and we claim that this implies that
$F \in B_p^{\gamma q}$ for some $\gamma \in \RR$.  To see this, let $\beta_{\nu}$ be as in the proof of 
Lemma \ref{besov2}; we need to show that 
$\{ 2^{\nu\gamma}\|\beta_{\nu}(\Delta)F\|_p: -1 \leq \nu \leq \infty\}$ is in $l^q$ for some $\gamma \in 
\RR$.  As in the proof of Lemma \ref{besov2}, for $\nu \geq -1$, let $J_{\nu}$ be the kernel of $\beta_{\nu}(\Delta)$. 
The arguments in the second paragraph of the proof of Lemma \ref{besov2} (specifically (\ref{jliest}), with $I = 0$,
in which case $J_{\nu}^I = J_{\nu}$), show that $\|J_{\nu}\|_{C^{M_0}({\bf M} \times {\bf M})} \leq C_02^{\nu(n+M_0)}$
for some $C_0 > 0$.
For any fixed $x \in {\bf M}$, let $J_{\nu,x}(y) = J_{\nu}(x,y)$.  Then, for any $x \in {\bf M}$,
\[ |[\beta_{\nu}(\Delta)F](x)| = |F(J_{\nu,x})| \leq C_12^{\nu(n+M_0)}, \]
for some $C_1 > 0$.  Thus 
\[ |[\beta_{\nu}(\Delta)F](x)|_p \leq C_1\mu({\bf M})^{1/p} 2^{\nu(n+M_0)}. \]
Therefore 
$\{ 2^{\nu\gamma}\|\beta_{\nu}(\Delta)F\|_p: -1 \leq \nu \leq \infty\}$ is in $l^q$ if $\gamma + n + M_0 < 0$,
that is, if $\gamma < -n -M_0$.  Fix $\gamma < \min(-n-M_0, \alpha)$; then $F \in B_p^{\gamma q}$,
and $\gamma \leq \alpha$.  In fact, since we are assuming that $F(1) = 0$, $F \in B_{p,0}^{\gamma q}$. \\ 

By Theorem \ref{besmain1}, we may choose $b_0$ sufficiently small that $S: B_{p,0}^{\alpha q}
\rightarrow B_{p,0}^{\alpha q}$ and $S: B_{p,0}^{\gamma q} \rightarrow B_{p,0}^{\gamma q}$ are both invertible if
$0 < b < b_0$.  For such $b$, since $S_0F \in B_{p,0}^{\alpha q} \subseteq B_{p,0}^{\gamma q}$, there is a
unique $F_1 \in B_{p,0}^{\alpha q}$ with $SF_1 = S_0F$, and a unique $F_2 \in B_{p,0}^{\gamma q}$ with
$SF_2 = S_0F$.  Now $S_0F$ is the sum (in $B_{p,0}^{\alpha q}$) of the series in (\ref{sumjopdfbes}).  Since
$F \in B_{p,0}^{\gamma q}$, that series converges in $B_{p,0}^{\gamma q}$ to $SF$.  Thus $S_0F = SF$
(as elements of $B_{p,0}^{\gamma q}$),
so $F_2 = F$.  Also 
$F_1 \in B_{p,0}^{\alpha q} \subseteq B_{p,0}^{\gamma q}$, and $SF_1 = S_0F$, so $F_1 = F_2 = F$.
But then $F = F_1 \in B_{p,0}^{\alpha q}$, as desired.  This completes the proof.

\begin{theorem}
\label{besmain3}
Say $c_0 , \delta_0, {\cal C} > 0$.  Fix an integer $l \geq 1$ with $2l > \max(n(1/p-1)_+ - \alpha, \alpha)$.
Say $f_0 \in {\cal S}(\RR^+)$, and let $f(s) = s^l f_0(s)$.  
Suppose also that the Daubechies condition (\ref{daub}) holds.
Then there exist constants $C_1 > 0$  and $0 < b_0 < 1$ as follows:\\

Say $0 < b < b_0$.  Then there exists a constant $C_2 > 0$ as follows:\\

Suppose that, for each $j$, we can write ${\bf M}$ as a finite disjoint union of measurable sets 
$\{E_{j,k}: 1 \leq k \leq N_j\}$, and that (\ref{diamleq}), (\ref{measgeq}) hold.  If $0 < p < 1$, we
assume that (\ref{measgeq2}) holds as well.
Select $x_{j,k} \in E_{j,k}$ for each $j,k$. 
  
For $t > 0$, let $K_t$ be the kernel of $f(t^2\Delta)$.  Set
\[ \varphi_{j,k}(y) =  \overline{K}_{a^j}(x_{j,k},y). \]
Then:\\
If $F \in B_{p,0}^{\alpha q}$, there exist constants $r_{j,k}$ with 
$(\sum_{j = -\infty}^{\infty} a^{-j\alpha q} [\sum_k \mu(E_{j,k})|r_{j,k}|^p]^{q/p})^{1/q} < \infty$
such that
\begin{equation}
\label{dfexp}
F = \sum_{j=-\infty}^{\infty} \sum_k \mu(E_{j,k}) r_{j,k} \varphi_{j,k},
\end{equation}
with convergence in $B_p^{\alpha q}$.   Further\\
\begin{equation}
\label{nrmeqiv2}
\|F\|_{B_p^{\alpha q}}/C_2 \leq 
\inf\{(\sum_{j = -\infty}^{\infty} a^{-j\alpha q} [\sum_k \mu(E_{j,k})|r_{j,k}|^p]^{q/p})^{1/q}:(\ref{dfexp}) \mbox{ holds}\}
 \leq C_1\|F\|_{B_p^{\alpha q}}.
\end{equation}
Moreover, if $p \geq 1$, then $C_2$ may be chosen to be independent of the choice of $b$
with $0 < b < b_0$.
\end{theorem}
{\bf Proof} We choose $b_0 > 0$ sufficiently small that $S: B_{p,0}^{\alpha q} \rightarrow B_{p,0}^{\alpha q}$
is invertible for $0 < b < b_0$.  For such $b$, if $F \in B_{p,0}^{\alpha q}$, then 
$S^{-1}F \in B_{p,0}^{\alpha q}$, so
\begin{equation}
\label{ssinvf}
F = S(S^{-1}F) = \sum_{j}\sum_k \mu(E_{j,k}) \langle  S^{-1}F, \varphi_{j,k}  \rangle  \varphi_{j,k}. 
\end{equation}
so that (\ref{dfexp}) holds with $r_{j,k} = \langle  S^{-1}F, \varphi_{j,k}  \rangle$.  By Theorem 
\ref{besmain2} and Theorem \ref{besmain1}, $(\sum_k \mu(E_{j,k})|r_{j,k}|^p]^{q/p})^{1/q} \leq 
C_1^{\prime} \|S^{-1}F\|_{B_p^{\alpha q}} \leq C_1\|F\|_{B_p^{\alpha q}}.$  

That leaves only the leftmost inequality in (\ref{nrmeqiv2}) to prove.  We need only show that, for any $F$
as in (\ref{dfexp}) (convergence in $B_p^{\alpha q}$), we have the inequality 
\begin{equation}
\label{frjkbpq}
\|F\|_{B_p^{\alpha q}} \leq C_2\sum_{j = -\infty}^{\infty} a^{-j\alpha q} [\sum_k \mu(E_{j,k})|r_{j,k}|^p]^{q/p})^{1/q}.
\end{equation}
with $C_2$ independent of $0 < b < b_0$ if $p \geq 1$.
It is enough to do this in each of two cases: (i) if $r_{j,k} = 0$ whenever $j > 0$; and (ii) if 
$r_{j,k} = 0$ whenever $j \leq 0$.  Case (i) follows at once from Lemma \ref{besov2} (and, if $0 < p < 1$, the hypotheses
(\ref{measgeq}), (\ref{measgeq2}), which imply (\ref{ejkro})).  Case (ii) follows from (*) of the proof of
Theorem \ref{besmain2}, since the inequality (\ref{cnmjpq}) there shows that $\|F\|_{C^N}$ is less than
or equal to the right side of (\ref{frjkbpq}), for any fixed $N$.  This completes the proof.

\ \\
STONY BROOK UNIVERSITY
 
\end{document}